\documentclass[12pt]{article} 
\usepackage{amsmath, amssymb, amsthm}
\usepackage{bbm}
\usepackage[normalem]{ulem}

\usepackage{latexsym}
\usepackage{graphicx}   
\usepackage{bm}
\usepackage{overpic}
\usepackage{harvard} 
  
\usepackage{exscale} 
\usepackage{amsfonts}
\usepackage{amssymb}  
\usepackage[usenames]{color} 

\definecolor{darkgreen}{rgb}{.3,.7,0}

 \def\ua{\uparrow}
 \def\da{\downarrow} 

 \def\wh{\widehat}
 \def\wt{\widetilde}

\def\bbar{\overline}
\def\bR{\mathbb{R}}

\def\cE{\mathcal E}  
\def\cC{\mathcal C}

\def\bR{\mathbb R}

\def\bE{\mathbb E}
\def\bN{\mathbb N}

\def\eps{\varepsilon}

\def\<{\langle}
\def\>{\rangle}

\def\Ind#1{\mathbbmss 1_{ #1 }} 

\def\ignore#1{{}}

\newcommand{\R}{\mathbb{R}}

\def\Proof{\bigskip\noindent{\bf Proof. }}
\newtheorem{theorem}{Theorem}

\newtheorem{lemma}{Lemma}
\newtheorem{corollary}{Corollary}

\theoremstyle{definition}

\newtheorem{example}{Example}
\newtheorem{remark}{Remark}

\begin{document}
\title{\LARGE\bf Capacitary measures for completely monotone kernels via singular control  }
\author{\hspace{-1cm} Aur\'elien Alfonsi\footnote{Universit\'e Paris-Est, CERMICS,
   Project team MathFi ENPC-INRIA-UMLV,
    Ecole des Ponts, 6-8 avenue Blaise Pascal,
   77455 Marne La Vall\'ee, France. {\tt alfonsi@cermics.enpc.fr
}}, {\setcounter{footnote}{6}}  Alexander
Schied\footnote{Department of Mathematics, University of Mannheim, 
A5, 6, 68131 Mannheim, Germany. 
{\tt schied@uni-mannheim.de}\hfill\break  The first author acknowledges the
support of the ``Chaire Risques Financiers'' of Fondation du Risque. The second author gratefully acknowledges support by Deutsche Forschungsgemeinschaft DFG} }
 
\date{\small January 11, 2012; this version: February 25, 2013
}
\maketitle

\abstract We give a singular control approach to the problem of minimizing an energy functional for measures with given total mass on a compact real interval, when  energy is defined in terms of a completely monotone kernel. This problem occurs both in potential theory and when looking for optimal financial order execution strategies under transient price impact. In our setup,  measures or order execution strategies are interpreted as singular controls, and the capacitary measure is the unique optimal control. The minimal energy, or equivalently the capacity of the underlying interval, is characterized by means of a nonstandard infinite-dimensional Riccati differential equation, which is analyzed in some detail. We then show that the capacitary measure has two Dirac components at the endpoints of the interval and a continuous Lebesgue density in between. This density can be obtained as the solution of a certain Volterra integral equation of the second kind.  

\bigskip


\medskip

\section{Introduction and statement of results}

\subsection{Background}

Let $G:\bR_+\to\bR_+$ be a function. The problem of minimizing the {energy functional}
$$\cE(\mu):=\frac12\int\int G(|t-s|)\,\mu(ds)\,\mu(dt)
$$
 over probability measures $\mu$ supported by a given compact set $K\subset\bR$ plays an important role in potential theory. A minimizing measure $\mu^*$, when it exists, is called a  {capacitary measure}, and the value $\text{Cap}\,(K):=1/\cE(\mu^*)$ is called the {capacity} of the set $K$; see, e.g., 
 \citeasnoun{choquet}, \citeasnoun{fuglede},   and \citeasnoun{landkof}.  See also \citeasnoun{AikawaEssen} or \citeasnoun{Helms} for more recent books on potential theory.  

In this paper, we develop a control approach to determining the capacitary distribution $\mu^*$ when $K$ is a compact interval and $G$ is a completely monotone function. In this approach,   measures $\mu$ on $K$ will be regarded as singular controls and $\cE(\mu)$ is the objective function. Our goal is to obtain qualitative structure theorems for the optimal control $\mu^*$ and characterize $\mu^*$  by means of certain differential and integral equations.

The intuition for  this control approach, and in fact   our original motivation, come from the problem of optimal order execution in mathematical finance.  In this problem, one considers an economic agent who wishes to liquidate a certain asset position of $x$ shares within the time interval $[0,T]$. This asset position can either be a long position ($x>0$) or a short position $(x<0$).  The order execution  strategy chosen by the investor is described by the asset position $X_t$ held  at time $t\in[0,T]$. In particular, one must have $X_0=x$.  Requiring the condition $X_{T+}=0$ assures that the initial position has been unwound by time $T$. The left-continuous path $X=(X_t)_{t\in[0,T]}$ will be nonincreasing for a pure sell strategy and nondecreasing for a pure buy strategy. A general strategy can consist of both buy and sell trades and hence can be described as the sum of a nonincreasing and a nondecreasing strategy. That is, $X$ is a path of finite variation. 

The problem the economic agent is facing is that his or her trades impact the price of the underlying asset. To model price impact, one starts by informally defining  $q\,dX_t$ as  the immediate price impact generated by the (possibly infinitesimal) trade $dX_t$ executed at time $t$. 
Next, it is an empirically well-established fact that price impact is transient and decays over time; see, e.g., \citeasnoun{MoroEtAl}. This decay of price impact can be described informally by requiring that $G(t-s)\,dX_s$ is the remaining impact at time $t$ of the impact generated by the trade $dX_s$. Here, $G:\bR_+\to\bR_+$ is a  nonincreasing function with $G(0)=q$, the {decay kernel}. Thus, $\int_{s<t}G(t-s)\,dX_s$ is the price impact of the strategy $X$, cumulated until time $t$. This price impact creates liquidation costs for the economic agent, and one can derive that, under the common  martingale assumption for unaffected asset prices, these costs are given by
\begin{equation}\label{costFunctional}
\cC(X):=\frac12\int_{[0,T]}\int_{[0,T]}G(|t-s|)\,dX_s\,dX_t
\end{equation}
plus a stochastic error term with expectation independent of the specific strategy $X$; see \citeasnoun{GSS}. 
Indeed, let us assume that asset prices are given by $S^X_t=S^0_t+ \int_0^tG(t-s)\,dX_s$ where $S^0$ is a continuous martingale and $\int_0^tG(t-s)\,dX_s$ models the price impact of the trading strategy at time~$t$. Then, we assume that the order $dX_t$ is made at the average price $\frac{1}{2}(S^X_{t-}+S^X_t)$ and costs $\frac{1}{2}(S^X_{t-}+S^X_t)\,dX_t$, which corresponds to a block shape limit order book, see~\citeasnoun{AFS2}. Accumulating these costs over $[0,T]$, integrating by parts twice, and taking expectations yields 
$$\mathbb{E} \left[\int_{[0,T]}\frac{1}{2}(S^X_{t-}+S^X_t)\,dX_t \right]=-S_0^0X_0+\bE[\,\cC(X)\,],
$$
where we have used  the fact that  $\mathbb{E} \left[\int_{[0,T]}S^0_tdX_t \right]=-S^0_0X_0$, due to the martingale assumption on $S^0$.
 Further details can be found in~\citeasnoun{GSS}.

Thus, minimizing the expected costs amounts to minimizing the functional $\cC(X)$ over all left-continuous strategies $X$ that are of bounded variation and satisfy $X_0=x$ and $X_{T+}=0$. This problem was  formulated and solved in the special case of exponential decay, $G(t)=e^{-\rho t}$, by \citeasnoun{ow}. The general case was analyzed by \citeasnoun{ASS} in discrete time and by \citeasnoun{GSS} in the continuous-time setup we have used above. We refer to \citeasnoun{AFS2},  \citeasnoun{AS}, \citeasnoun{GSS2}, \citeasnoun{PredoiuShaikhetShreve}, \citeasnoun{SchiedSlynko}, and \citeasnoun{GatheralSchiedSurvey} for further discussions and additional references in the context of mathematical finance.

Clearly, the cost functional $\cC(X)$ coincides with the energy functional $\cE(\nu^X)$ of the measure $\nu^X(dt):=d X_t$. So finding an optimal order execution strategy is basically equivalent to determining a capacitary measure for $[0,T]$.  There is one important difference, however: capacitary measures are determined as minimizers of $\cE(\mu)$ with respect to all \emph{nonnegative} measures $\mu$ on $[0,T]$ with total mass $1$, while $\nu^X$ may be a \emph{signed} measure with given total mass $\nu^X([0,T])=-x$. This difference can become significant if $G(|\cdot|)$ is only required to be positive definite in the sense of Bochner (which is essentially equivalent to $\cC(X)\ge0$ for all $X$), because then minimizers of the unconstrained problem need not exist.  It was first shown by \citeasnoun{ASS}, and later extended to continuous time by \citeasnoun{GSS},
that  a unique optimal order execution strategy $X^*$ exists and that $X^*$ is a monotone function of $t$ when $G$ is convex and nonincreasing.  This result has the important consequence that the constrained problem of finding a capacitary measure is equivalent to the unconstrained problem of  determining an optimal order execution strategy. 

 In this paper, we aim at describing the structure of capacitary measures/optimal order execution strategies. To this end, it is instructive to first look at two specific examples in which the optimizer is known in explicit form. \citeasnoun{ow} find that for exponential decay, $G(t)=e^{-\rho t}$, the capacitary measure $\mu^*$ has two singular components at $t=0$ and $t=T$ and a constant Lebesgue density on $(0,T)$:
\begin{equation}\label{strategy_OW}
\mu^*(dt)=\frac1{2+\rho T}\,\delta_0(dt)+\frac\rho{2+\rho T}\,dt+\frac1{2+\rho T}\,\delta_T(dt).
\end{equation}
Numerical experiments show that it is a common pattern that capacitary measures for nonincreasing convex kernels have two singular components at $t=0$ and $T=0$ and a Lebesgue density on $(0,T)$.  However, the capacitary measure for $G(t)=\max\{0,1-\rho t\}$ is the purely discrete measure
$$\mu^*=\frac{1}{2+N}\sum_{i=0}^N\Big(1-\frac i{N+1}\Big)\big(\delta_{\frac i\rho}+\delta_{T-\frac i\rho}\big),$$
where $N:=\lfloor \rho T\rfloor$ \cite[Proposition 2.14]{GSS}. 

So it is an interesting question for which nonincreasing, convex kernels $G$ the capacitary measure $\mu^*$  has singular components only at $t=0$ and $t=T$ and is (absolutely) continuous on $(0,T)$.  It turns out that a sufficient condition is the \emph{complete monotonicity} of $G$, i.e., $G$ belongs to $C^\infty((0,\infty))$ and $(-1)^nG^{(n)}$ is nonnegative in $(0,\infty)$ for $n\in\bN$. More precisely, we have the following result, which is in fact an immediate corollary of the main results in this paper.

\begin{corollary}Suppose that $G$ is completely monotone  with
  $G''(0+):=\lim_{t\da0}G''(t)<\infty$. Then the capacitary measure $\mu^*$
  has two Dirac components at $t=0$ and $t=T$ and is has a continuous Lebesgue density on $(0,T)$.
\end{corollary}

\subsection{Statement of main results}

Our main results do not only give the preceding qualitative statement on the form of $\mu^*$ but they also provide quantitative descriptions of the Dirac components of $\mu^*$ and of its Lebesgue density on $(0,T)$.  To prepare for the statement of these results, let us first assume that $G(0)=1$, which we can do without loss of generality. Then we recall that by the celebrated Hausdorff--Bernstein--Widder theorem \cite[Theorem IV.12a]{Widder}, $G$ is completely monotone if and only if it is the Laplace transform of a Borel probability measure $\lambda$ on $\bR_+$:
$$G(t)=\int e^{-\rho t}\,\lambda(d\rho), \qquad t\ge0. 
$$
In particular, every exponential polynomial,
\begin{equation}\label{pol_G}
G(t)=\sum_{i=0}^d\lambda_ie^{-\rho_i t},
\end{equation}
with $\lambda_i,\,\rho_i\ge0$ and $\sum_i\lambda_i=1$ is completely monotone. Another example is power-law decay,
$$G(t)=\frac{1}{(1+t)^\gamma}\qquad\text{for some $\gamma>0$,}
$$
 which is a popular choice for the decay of price impact in the econophysics literature; see \citeasnoun{Gatheral} and the references therein. 
We assume henceforth that $G''(0+)<\infty$, which is equivalent to  
\begin{equation}\label{bbarRho2Cond}
\bbar\rho:=\int\rho\,\lambda(d\rho)<\infty\qquad\text{and}\qquad \bbar{\rho^2}:=\int\rho^2\,\lambda(d\rho)<\infty.
\end{equation}

\bigskip

A crucial role will be played by the following  infinite-dimensional Riccati equation for functions $\varphi:[0,\infty)\times\bR_+^2\to\bR$,
\begin{equation}\label{InfiniteDimRiccatiEqn}
\varphi'(t,\rho_1,\rho_2)+(\rho_1+\rho_2)\varphi(t,\rho_1,\rho_2)=\frac1{2\bbar\rho}\Big(\rho_1+\int x\varphi(t,\rho_1,x)\,\lambda(dx)\Big)\Big(\rho_2+\int x\varphi(t,x,\rho_2)\,\lambda(dx)\Big)
\end{equation} 
where $\varphi'$ denotes the time derivative of $\varphi$, and the function $\varphi$ satisfies the initial condition
\begin{equation}\label{InfiniteDimRiccatiEqnInitialCond}
\varphi(0,\rho_1,\rho_2)=1\qquad\text{for all $\rho_1,\rho_2 \ge 0$.}
\end{equation} 

\begin{remark}When writing \eqref{InfiniteDimRiccatiEqn} in the form $\varphi'=F(\varphi)$ one sees that the functional $F$ is {not a continuous map}  from some reasonable function space into itself, unless $\lambda$ is concentrated on a compact interval. For instance, it involves the typically unbounded linear operator $\varphi\mapsto(\rho_1+\rho_2)\varphi$. Therefore, existence and uniqueness of solutions to~\eqref{InfiniteDimRiccatiEqn},~\eqref{InfiniteDimRiccatiEqnInitialCond} does not follow by an immediate application of standard results such as the Cauchy--Lipschitz/Picard--Lindel\"of theorem in Banach spaces  \cite[Theorem 3.4.1]{HillePhillips} or more recent ones such as those in \citeasnoun{Teixeira} and the references therein. In fact, even in the simplest case in which $\lambda$ reduces to a Dirac measure, the existence of global solution hinges on the initial condition; it is easy to see that solutions blow up when $\varphi(0)$ is not  chosen in a suitable manner.
\end{remark}

We now state a  result on the global existence and uniqueness of \eqref{InfiniteDimRiccatiEqn}, \eqref{InfiniteDimRiccatiEqnInitialCond}. It states that the solution takes values in the locally convex space  $C(\bR^2_+)$ endowed with topology of locally uniform convergence. For integers $k\ge0$, the space $C^k([0,\infty); C(\bR^2_+))$ will consist of all continuous functions $\varphi:[0,\infty)\to C(\bR^2_+)$ which, when considered as functions $\varphi:[0,\infty)\to C(K)$ for some compact subset $K$ of $\bR^2_+$, belong to $C^k([0,\infty); C(K))$.

\begin{theorem}\label{ODEthm}When $G''(0+)<\infty$ the initial value problem~\eqref{InfiniteDimRiccatiEqn},~\eqref{InfiniteDimRiccatiEqnInitialCond} admits a unique solution $\varphi$ in the class of functions $\wt\varphi$ in $ C^1([0,\infty);C(\bR^2_+))$ that satisfy an inequality of the form
\begin{equation}\label{phiiEstimateEq}
0\le\wt\varphi(t,\rho_1,\rho_2)\le  c(1+\rho_1)(1+\rho_2),
\end{equation}
where $c$ is a constant that may depend on $\varphi$ and locally uniformly on $t$. 
Moreover,  $\varphi$ has the following properties.
\begin{enumerate}
\item $\varphi$ is strictly positive.

\item\label{SymmetryProperty} $\varphi$ is symmetric: $\varphi(t,\rho_1,\rho_2)=\varphi(t,\rho_2,\rho_1)$ for all $(\rho_1,\rho_2)\in\bR_+^2$.
\item $1=\int\varphi(t,\rho,x)\,\lambda(dx)=\int\varphi(t,x,\rho)\,\lambda(dx)$ for all $\rho\ge0$.
\item\label{CinftyProperty}$\varphi\in C^{2}([0,\infty);C(\bR^2_+))$.
\item For every $t$, the kernel $\varphi(t,\cdot,\cdot)$  is nonnegative definite on $ L^2(\lambda)$, i.e., 
\begin{equation}\label{nonneg}
\int\int f(x)f(y)\varphi (t,x,y)\,\lambda(dx)\,\lambda(dy)\ge0\qquad\text{for $f\in L^2(\lambda)$.}
\end{equation}
\item The functions $\varphi(t,\rho_1,\rho_2)$ and $\varphi'(t,\rho_1,\rho_2)$ satisfy local Lipschitz conditions in $(\rho_1,\rho_2)$, locally uniformly in $t$.
\end{enumerate}
\end{theorem}

\bigskip

In Section \ref{sec_computational} we will discuss computational aspects of the initial value problem~\eqref{InfiniteDimRiccatiEqn},~\eqref{InfiniteDimRiccatiEqnInitialCond}. In particular, we will discuss its solution when $G$ is an exponential polynomial of the form \eqref{pol_G} and we will provide closed-form solutions in the cases $d=1$ and $d=2$.

We can now explain how to use singular control in approaching the minimization
of $\cE(\mu)$ or $\cC(X)$. To this end, using  order execution strategies
$X=(X_t)$ will be more convenient than using the formalism of the associated measures $\mu(dt)=dX_t$ because of the natural dynamic interpretation of $X$. Henceforth, a $[0,T]$-admissible strategy  will be a left-continuous function $(X_t)$ of bounded variation such that $X_{T+}=0$. Our goal is to minimize the cost functional $\cC(X)$ defined in \eqref{costFunctional} over all $[0,T]$-admissible strategies with fixed initial value $X_0=x$. 
Clearly, this problem is not yet  suitable for the application of control techniques since $\cC(X)$ depends on the entire path of $X$. We therefore introduce the auxiliary functions
\begin{equation}\label{def_EX}
E^X_t(\rho):=\int_{[0,t)}e^{-\rho(t-s)}\,dX_s,\qquad \text{for }\rho\ge0.
\end{equation} 
These functions will play the role of state variables that are controlled by the strategy $X$. 

\begin{lemma}\label{CostFunctionalLemma}For any $[0,T]$-admissible strategy $X$, the function $E_t^X(\rho)$ is uniformly bounded in $\rho$ and $t$. Moreover,
\begin{equation}\label{cCEFormulaEq}
\cC(X)=\int_{[0,T)}\int E_t^X(\rho)\,\lambda(d\rho)\,dX_t+\frac12\sum_{t\le T}(\Delta X_t)^2,
\end{equation} 
where $\Delta X_t:=X_{t+}-X_t$ denotes the jump of $X$ at $t$.
\end{lemma}

\Proof Clearly, $|E_t^X(\rho)|\le\|X\|_{\text{var}}$, where $\|X\|_{\text{var}}$ denotes the total variation of $X$ over
  $[0,T]$.  To obtain \eqref{cCEFormulaEq}, we integrate by parts to get
$$\cC(X)=\int_{{[0,T)}}\int_{{[0,t)}}{G}(t-s)\,dX_s\,dX_t+\frac{{G}(0)}2\sum_{t\le
  T}(\Delta X_t)^2.
$$
Now we write $G(t-s)$ as $\int e^{-\rho(t-s)}\,\lambda(d\rho)$ and apply Fubini's theorem. \qed.

\bigskip

The form \eqref{cCEFormulaEq} of our cost functional is now suitable for the application of control techniques.  To state  our main result, we let $\varphi$ be the solution of our infinite-dimensional Riccati equation as provided by Theorem \ref{ODEthm} and we define
\begin{equation}\label{psiDefinitionEq}
\varphi_0(t):=\varphi(t,0,0)\qquad\text{and}\qquad \psi(t,\rho):=\int x\varphi(t,x,\rho)\,\lambda(dx)
\end{equation}

\bigskip

\begin{theorem}\label{verificationThm}Let $X^*$ be the unique optimal strategy in the class of $[0,T]$-admissible strategies with initial value $X_0=x$. Then
\begin{equation}\label{MinimalCostEq}
\cC(X^*)=\frac{x^2}{2\varphi_0(T)}.
\end{equation}
Moreover, $X^*$  has jumps at  $t=0$ and $t=T$    of size
$$\Delta X_0^*=\Delta X_T^*=-\frac{\psi(T,0)}{2\bbar\rho\varphi_0(T)}\,x
$$
and is continuously differentiable on $(0,T)$. The derivative $\xi(t)=\frac d{dt}X^*_t$ is the unique continuous solution of the  Volterra integral equation 
\begin{equation}\label{VolterraEqn}
\xi (t)=f(t)+\int_0^t K(t,s)\xi (s)\,ds,
\end{equation}
where, for 
\begin{eqnarray} \Theta(t,\rho):=\frac{\rho+\psi(t,\rho)}{\psi(t,0)}\int x^2\varphi(t,x,0)\,\lambda(dx)-\int x^2\varphi(t,x,\rho)\,\lambda(dx)+\rho^2,\label{ThetaEq}
\end{eqnarray}
the function $f$ and the kernel $K(\cdot,\cdot)$ are given by
\begin{equation}
f(t)=\frac{\Delta X_0^*}{2\bbar\rho}\int e^{-\rho t}\Theta(T-t,\rho)\,\lambda(d\rho),\qquad K(t,s)=\frac1{2\bbar\rho}\int e^{-\rho (t-s)} \Theta(T-t,\rho)\,\lambda(d\rho).
\end{equation}
\end{theorem}

\bigskip

Let us recall that we know in addition from Theorem~2.20 in~\citeasnoun{GSS} that
$t\in[0,T] \mapsto X^*_t$ is monotone. 
The identity \eqref{MinimalCostEq} immediately yields the following formula for the capacity of a compact interval.

\begin{corollary}If $G''(0+)<\infty$, the capacity of a compact interval $[a,b]$ is given by
$$\text{\rm Cap}\,([a,b])=2\varphi_0(b-a).
$$
\end{corollary}

\subsection{Computational aspects}\label{sec_computational}

In general, the Riccati
equation~\eqref{InfiniteDimRiccatiEqn},~\eqref{InfiniteDimRiccatiEqnInitialCond} cannot be solved explicitly. 
A closed-form solution exists, however, when $G$ is an exponential polynomial as in \eqref{pol_G}, i.e.,  when $\lambda$ has a discrete support. Let
us assume that $\lambda(dx)=\sum_{i=0}^d \lambda_i \delta_{\rho_i}(dx)$, with
$\rho_0=0<\rho_1<\dots<\rho_d$, $\lambda_i\ge 0$, and $\sum_{i=0}^d\lambda_i=1$. All the input that is needed in Theorem \ref{verificationThm} are the values 
$\varphi_{ij}(t):=\varphi(t,\rho_i,\rho_j)$, for $0\le i,j\le d$.
By Theorem \ref{ODEthm},  $\varphi(t)$ is a symmetric matrix  that solves the
following  matrix Riccati equation:
\begin{equation}\label{Matrix Ricatti eqn}
\varphi'=-\varphi M^{(3)}\varphi-\varphi M^{(4)}+M^{(1)}\varphi+M^{(2 )},
\end{equation}
with $M^{(3)}_{ij}=-\frac1{{2\bbar\rho}}\lambda_i\rho_i\lambda_j\rho_j
$, $M^{(4)}_{ij}=-\frac{\lambda_i\rho_i\rho_j}{2 \bbar\rho}+\delta_{ij}\rho_i
$, $M^{(1)}=-(M^{(4)})^\top$ and $M^{(2)}_{ij}=\frac{\rho_i\rho_j}{2 \bbar\rho}
$. According to \citeasnoun{Levin}, the solution of this equation is given by
$$\varphi(t)=(R^{(1)}(t)\mathbf{1}+R^{(2)}(t))(R^{(3)}(t)  \mathbf{1}+R^{(4)}(t))^{-1},$$
where $\mathbf{1}_{ij}=1$  and 
$$ R(t)=\left[\begin{matrix}
    R^{(1)}(t) & R^{(2)}(t) \\R^{(3)}(t) & R^{(4)}(t)
  \end{matrix}
\right]= \exp\left(t \left[\begin{matrix}
    M^{(1)} & M^{(2)} \\ M^{(3)} & M^{(4)}
  \end{matrix} \right]\right).$$
In the special cases $d=1$ and $d=2$, the solution of the Riccati
equation~\eqref{InfiniteDimRiccatiEqn},~\eqref{InfiniteDimRiccatiEqnInitialCond} 
becomes even easier and, to some extend, becomes explicit. We demonstrate this first for $d=1$ and then for $d=2$:

\begin{example} In the case $d=1$, $G$ is of the form $G(t)=\lambda+(1-\lambda)e^{-\rho t}$ for some $\lambda\in[0,1)$ and some $\rho>0$. Clearly, we can set $\lambda:=0$ without changing the optimization problem.  Then $\bbar\rho=\rho_1=\rho$,
and \eqref{InfiniteDimRiccatiEqn} becomes
\begin{eqnarray*}
\varphi'_{00}=\frac\rho2\,\varphi_{01}^2,\quad
\varphi'_{01}+\rho\varphi_{01}=\frac\rho2\, (1+\varphi_{11})\varphi_{01},\quad
\varphi'_{11}+2\rho\varphi_{11}=\frac\rho2\,(1+\varphi_{11})^2.
\end{eqnarray*}
For the initial condition $\varphi_{kl}(0)=1$, the preceding equation has the unique solution $\varphi_{11}\equiv\varphi_{01}\equiv 1$ and $\varphi_{00}(t)=1+\rho t/2$. The condition~\eqref{X*Cond2} thus reduces to
$0={X_t}+E^X_1(t)\big(1+\rho(T-t)\big)$, 
which easily yields \eqref{strategy_OW} as unique solution.
\end{example} 

\begin{example}In the case $d=2$, we can assume that $G$ is of the form $G(t)=\lambda_1e^{-\rho_1t}+\lambda_2e^{-\rho_2t}$, where $\lambda_1+\lambda_2=1$.
Consider a solution $\varphi_{ij}$ $(i,j=0,\dots,2)$ of the matrix Riccati equation \eqref{Matrix Ricatti eqn} with $\lambda_0=0$.  We can simplify \eqref{Matrix Ricatti eqn} by using the relation
\begin{eqnarray}
\lambda_1\varphi_{i1}+\lambda_2\varphi_{i2}=1,\qquad i=0,\dots, 2.\label{d=2SimplificationEq}
\end{eqnarray}
Indeed, the equation for $\varphi_{11}$ then becomes
\begin{eqnarray*}
\varphi_{11}'+2\rho_1\varphi_{11}=\frac1{2\bbar\rho}\big(\rho_1+\lambda_1\rho_1\varphi_{11}+\lambda_2\rho_2\varphi_{12}\big)^2=\frac1{2\bbar\rho}\big(\rho_1+\rho_2+\lambda_1(\rho_1-\rho_2)\varphi_{11}\big)^2.
\end{eqnarray*}
This is an autonomous ODE that, for the initial condition $\varphi_{11}(0)=1$,  is solved by
\begin{equation}\label{phi11}
\varphi_{11}(t)=c_1+\bigg[\Big(\frac1{1-c_1}-\frac14\frac{\lambda_1^2(\rho_1-\rho_2)^2}{\sqrt{\rho_1\rho_2\bbar\rho\check\rho}}\Big)\exp\Big(\frac{2\sqrt{\rho_1\rho_2\bbar\rho\check\rho}}{\bbar\rho}\cdot t\Big)+\frac14\frac{\lambda_1^2(\rho_1-\rho_2)^2}{\sqrt{\rho_1\rho_2\bbar\rho\check\rho}}\bigg]^{-1},
\end{equation}
where
$$\check\rho:=\lambda_1\rho_2+\lambda_2\rho_1 \text{ and } c_1=\Big(\frac{\sqrt{\rho_1\bbar\rho}-\sqrt{\rho_2\check\rho}}{\lambda_1(\rho_1-\rho_2)}\Big)^2.
$$
We can notice that $\varphi_{11}(+\infty)=c_1$ and $\varphi_{11}(-\infty)=\Big(\frac{\sqrt{\rho_1\bbar\rho}+\sqrt{\rho_2\check\rho}}{\lambda_1(\rho_1-\rho_2)}\Big)^2$.

Similarly, 
$$\varphi_{22}'+2\rho_2\varphi_{22}=\frac1{2\bbar\rho}\big(\rho_1+\rho_2+\lambda_2(\rho_2-\rho_1)\varphi_{22}\big)^2,
$$
which for the initial condition $\varphi_{22}(0)=1$  is solved by
\begin{equation}\label{phi22}
\varphi_{22}(t)=c_2+\bigg[\Big(\frac1{1-c_2}-\frac14\frac{\lambda_2^2(\rho_1-\rho_2)^2}{\sqrt{\rho_1\rho_2\bbar\rho\check\rho}}\Big)\exp\Big(\frac{2\sqrt{\rho_1\rho_2\bbar\rho\check\rho}}{\bbar\rho}\cdot t\Big)+\frac14\frac{\lambda_2^2(\rho_1-\rho_2)^2}{\sqrt{\rho_1\rho_2\bbar\rho\check\rho}}\bigg]^{-1},
\end{equation}
where
$$c_2=\Big(\frac{\sqrt{\rho_2\bbar\rho}-\sqrt{\rho_1\check\rho}}{\lambda_2(\rho_1-\rho_2)}\Big)^2.
$$
From \eqref{d=2SimplificationEq} we can now easily compute $\varphi_{12}$. 

Next, using once again \eqref{d=2SimplificationEq}, we find that  $\varphi_{01}$ solves
$$\varphi_{01}'+\rho_1\varphi_{01}=\frac1{2\bbar\rho}\big(\rho_2+\lambda_1(\rho_1-\rho_2)\varphi_{01}\big)\big(\rho_1+\rho_2+\lambda_1(\rho_1-\rho_2)\varphi_{11}\big).
$$
That is,
\begin{equation}\label{ODEphi01}
\varphi_{01}'+\left[\frac{\rho_1 \bbar\rho+ \rho_2\check\rho}{2 \bbar\rho}-\frac{\lambda_1^2}{2 \bbar\rho}(\rho_1-\rho_2)^2 \varphi_{11} \right]\varphi_{01}=\frac{\rho_2}{2 \bbar\rho}\left[\rho_1+\rho_2+\lambda_1(\rho_1-\rho_2)\varphi_{11} \right].
\end{equation}
We set $B_1=\frac14\frac{\lambda_1^2(\rho_1-\rho_2)^2}{\sqrt{\rho_1\rho_2\bbar\rho\check\rho}}$, $A_1=\frac1{1-c_1}-B_1$,  and $k=\frac{\sqrt{\rho_1\rho_2\bbar\rho\check\rho}}{\bbar\rho}$, so that $\varphi_{11}(t)=c_1+\frac{1}{A_1e^{2kt}+B_1}$. Then, we can check that $\frac{1}{A_1e^{kt}+B_1e^{-kt}}$ is a solution of the fundamental system. By using a variation of parameters, we get that the solution of~\eqref{ODEphi01} satisfying  $\varphi_{01}(0)=1$ is given by
$$ \varphi_{01}(t)=\frac{A_1\varphi_{01}(+\infty)e^{kt}+B_1\varphi_{01}(-\infty)e^{-kt}+C_{01}}{A_1e^{kt}+B_1e^{-kt}},$$
with $$ \varphi_{01}(\pm \infty)=\frac{\pm \sqrt{\rho_1\rho_2\bbar\rho\check\rho}-\rho_2\check\rho }{\lambda_1 \check\rho (\rho_1-\rho_2)} \text{ and } C_{01}=A_1(1-\varphi_{01}(+\infty)) +B_1(1-\varphi_{01}(-\infty)).$$
 Then, $\varphi_{02}$ can be easily deduced from~\eqref{d=2SimplificationEq}.

It remains  to compute $\varphi_{00}$, which solves
$$  \varphi_{00}'=\frac{1}{2 \bbar\rho}\left[ \rho_2+\lambda_1(\rho_1-\rho_2)\varphi_{01}\right]^2. $$
We set $\tilde{C}_{01}=\lambda_1(\rho_1-\rho_2)C_{01}$ and get after some calculations:
\begin{eqnarray*}  
\varphi_{00}'(t)&=&\frac{1}{2 \bbar\rho}\left[ \frac{A_1k\frac{\bbar \rho}{\check \rho}e^{kt}-B_1k\frac{\bbar \rho}{\check \rho}e^{-kt} +\tilde{C}_{01}}{A_1e^{kt}+B_1e^{-kt}}  \right]^2 \\
&=&\frac{\rho_1 \rho_2}{2 \check \rho} + \frac{\tilde{C}_{01}}{\check \rho} \frac{A_1ke^{kt}-B_1ke^{-kt}}{(A_1e^{kt}+B_1e^{-kt})^2}+\frac{\tilde{C}_{01}^2\check \rho-4A_1B_1\rho_1 \rho_2 \bbar \rho}{2 \bbar \rho \check \rho}\frac{1}{(A_1e^{kt}+B_1e^{-kt})^2}.
\end{eqnarray*} 
Thus, we finally get:
\begin{eqnarray}\label{ODEphi00}
\varphi_{00}(t)&=&1+\frac{\rho_1 \rho_2}{2 \check \rho}t - \frac{\tilde{C}_{01}}{\check \rho}\left(\frac{1}{A_1e^{kt}+B_1e^{-kt}} -\frac{1}{A_1+B_1} \right) \\
&&+\frac{\tilde{C}_{01}^2\check \rho-4A_1B_1\rho_1 \rho_2 \bbar \rho}{4 B_1 k \bbar \rho \check \rho}\left(\frac{e^{kt}}{A_1e^{kt}+B_1e^{-kt}} -\frac{1}{A_1+B_1} \right). \nonumber
\end{eqnarray}
This completes this example. 
\end{example}

\bigskip

Given the solution $\varphi$ of the Riccati equation, we can approximate
the continuous time strategy by a discrete one as follows ($x_i$ will denote
the trading size at time $iT/N$).
\begin{itemize}
\item  We first set $x_0=\frac{\psi(T,0)}{2\bbar\rho\varphi_0(T)}$ and
  $E_0(\rho_\ell)=x_0$, $0\le \ell \le d$.
\item Suppose that $1\le i< N$ and that $x_{i-1}$ and $E_{i-1}(\rho_\ell)$ have been computed. Then, we
  set thanks to~\eqref{X*Cond2}:
  $$x_i=1-\sum_{j=0}^{i-1}x_j-\int E_{i-1}(\rho)e^{-\rho
    T/N}\theta(T-iT/N,\rho)\lambda(d\rho), \ E_i(\rho_\ell)= E_{i-1}(\rho_\ell)e^{-\rho_\ell
    T/N}+x_i.$$
\item Set $x_N=1-\sum_{j=0}^{i-1}x_j$.
\end{itemize}

Alternatively, we could have approximated  the minimization of the
cost~\eqref{costFunctional} by the following discrete problem. Let  $M_{i,j}=G\left(|i-j|\frac{T}{N}\right)$, $0\le i,j\le N$, and consider
\begin{equation}\label{Algo1}\text{ minimize $\frac{1}{2} x^T M x$ over $x \in \R^{N+1}$ s.t
  $\sum_{i=0}^N x_i=1$}.
\end{equation}
The solution of this problem is obviously given by
$\frac{1}{\mathbf{1}^TM^{-1} \mathbf{1}}M^{-1} \mathbf{1}$, where
$\mathbf{1}_i=1$ for $0\le i\le N$.  From a financial point of view, the minimization problem~\eqref{Algo1} gives the optimal strategy when it is only possible to trade at the times~$iT/N$, while the original problem~\eqref{costFunctional} allows to trade continuously.         In potential theory, it corresponds to computing the capacitary distribution of the set $\{iT/N\,|\,i=0,\dots, N\}$. It was shown in the proof of Theorem 2.20 in \citeasnoun{GSS} that for   $N\uparrow\infty$ these cap2acitary distributions converge in the weak topology of probability measures  to the capacitary distribution $dX^*$ constructed in~Theorem~\ref{verificationThm}. Explicit solutions of \eqref{Algo1} for the choices $G(t)=e^{-\rho t}$ and $G(t)=(1-\rho t)^+$ were given in \citeasnoun{AFS1} and \citeasnoun{ASS} (note, however, that $G(t)=(1-\rho t)^+$ is not completely monotone).

\begin{figure}[htbp]
 \centering
 \includegraphics[width=0.8\linewidth]{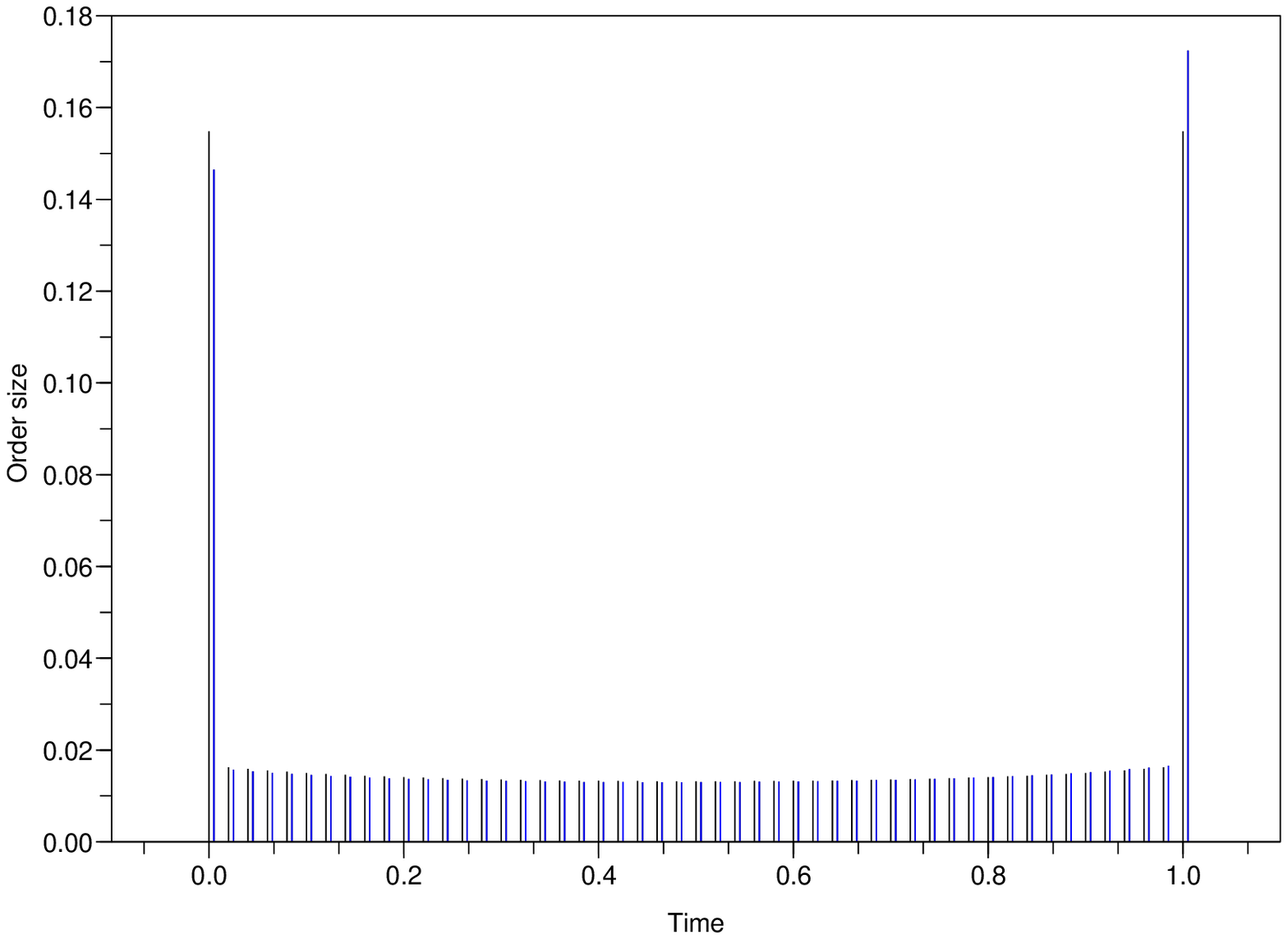} 
 \caption{Comparison of the approximated optimal strategies $(x_i,0\le i \le N)$ obtained
   with~\eqref{Algo1} and with the  method based on the Riccati equation (slightly shifted to the right) .}
 \label{figure_strategies}

 \centering
 \includegraphics[width=0.8\linewidth]{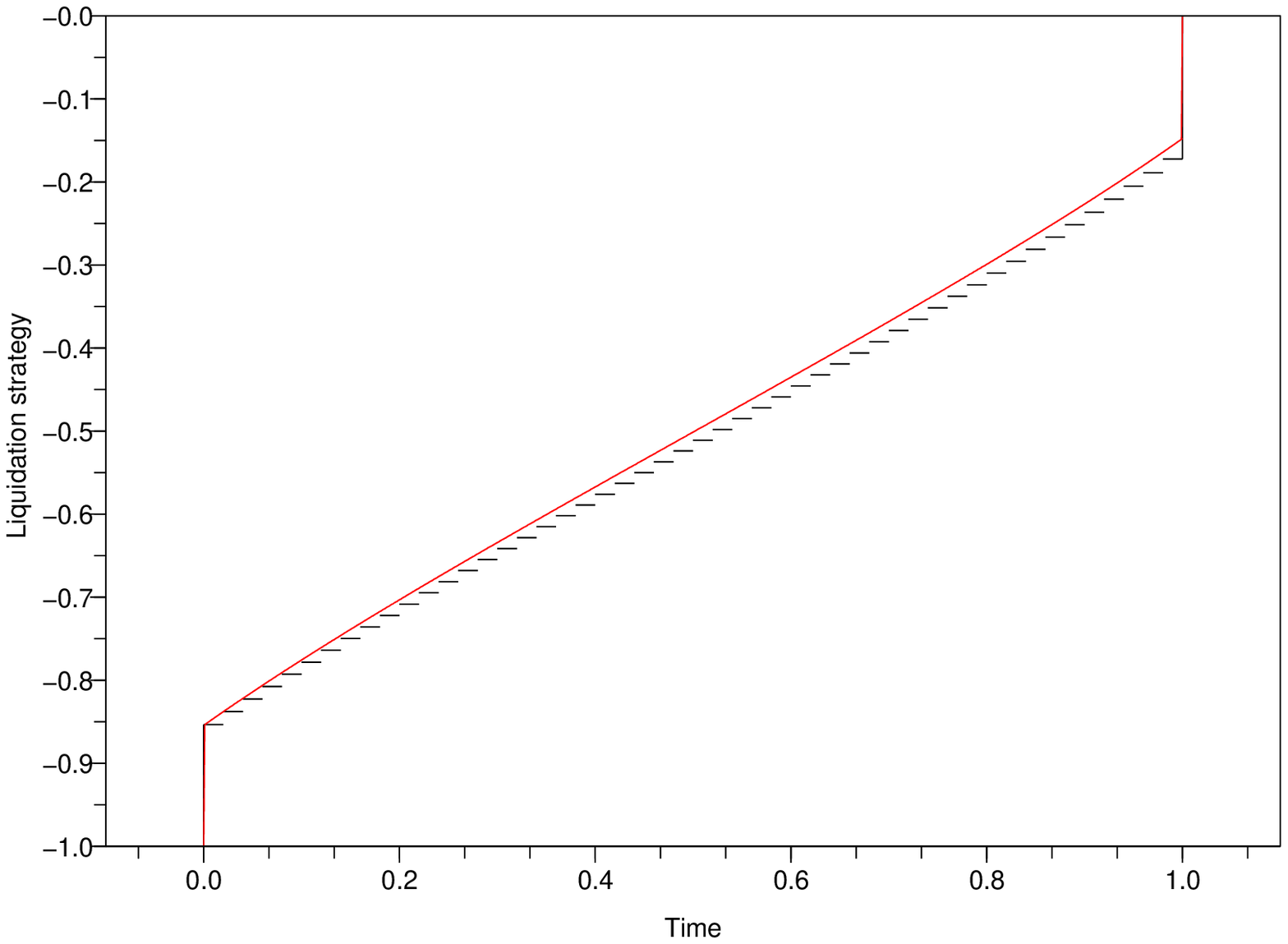} 
 \caption{Comparison of the approximated optimal strategies given by
   the Riccati method and $N=50$, 
   with the optimal continuous one~$X^*_t$ (computed with $N=1000$) .}
 \label{figure_strategies2}
\end{figure}

We have computed and plotted the solutions given by both methods in
Figure~\ref{figure_strategies} for $T=1$, $N=50$, and
$\lambda(d\rho)=0.1\delta_0(d\rho)+0.2\delta_1(d\rho)+0.2\delta_3(d\rho)+0.2\delta_5(d\rho)+0.2\delta_7(d\rho)+0.1\delta_{10}(d\rho)$. They
are already rather close together for $N=50$, and they merge when $N\rightarrow
+\infty$. Let us discuss briefly the time complexity of the two methods. The
 one given by~\eqref{Algo1} gets very slow when $N$ gets
large since it involves the inversion of a $N\times N$ matrix. Instead, when
$\lambda$ has a discrete support, the matrix Riccati equation can be solved
quickly and the algorithm above has a $O(N)$ time complexity, which is much
faster. However, this is no longer true when $\lambda$ does not have discrete
support. In that case, we have to approximate $\lambda$ by a discrete measure,
which means that we have to  increase~$d$. Doing so, will slow down the  algorithm based on the Riccati equation. A rigorous treatment of the
convergence rate and time complexity of both algorithms is beyond the scope of
this paper and is left for future research. 

\section{Proofs}

\subsection{Proof of Theorem \ref{ODEthm}}

Let us write~\eqref{InfiniteDimRiccatiEqn} in the form $\varphi'(t)=F_\lambda(\varphi(t))$, where
\begin{equation}\label{F(f)Eq}
F_\lambda(f)(\rho_1,\rho_2)=-(\rho_1+\rho_2)f(\rho_1,\rho_2)+\frac1{2\bbar\rho}\Big(\rho_1+\int x f(\rho_1,x)\,\lambda(dx)\Big)\Big(\rho_2+\int xf(x,\rho_2)\,\lambda(dx)\Big).
\end{equation}

\begin{lemma}\label{RiccatiLemma1}Suppose that $\lambda$ is supported by the compact interval $[0,\rho_{\max}]$. Then~\eqref{InfiniteDimRiccatiEqn},~\eqref{InfiniteDimRiccatiEqnInitialCond} admits a unique solution $\varphi\in C^1([0,\infty); C(\bR^2_+))$. Moreover, $\varphi$  has the properties {\rm (a)}, {\rm (b)}, and {\rm (c)} in the statement of Theorem~\ref{ODEthm}.
\end{lemma}

\Proof  Let $J\subset\bR_+$ be any compact interval containing $[0,\rho_{\max}]$. 
Then $F_\lambda$ defined in~\eqref{F(f)Eq} maps  $C(J\times J)$ into itself. Moreover, $F_\lambda$  is  Lipschitz continuous with respect to the sup-norm on every bounded subset of $C(J\times J)$. Hence, the Cauchy--Lipschitz/Picard--Lindel\"of theorem in Banach spaces implies the existence of a unique local solution $\varphi_J\in C^1([0,t_J); C(J\times J))$ for some maximal time $t_J>0$ \cite[Theorem 3.4.1]{HillePhillips}. 
We will show below  that $t_J=\infty$.  
Then, if $\wt J\supset J$ is another compact interval, the restriction of $\varphi_{\wt J}(t)$ to $J$ must coincide with  $\varphi_{ J}(t)$ due to the uniqueness of solutions. This consistency then implies the existence and uniqueness of solutions $\varphi\in C^1([0,\infty); C(\bR^2_+)))$.
Moreover, the uniqueness of solutions and the fact that both~\eqref{InfiniteDimRiccatiEqn} and~\eqref{InfiniteDimRiccatiEqnInitialCond}  are symmetric in $\rho_1$ and $\rho_2$ implies that  $\varphi(t,\rho_1,\rho_2)=\varphi(t,\rho_2,\rho_1)$ for all $(\rho_1,\rho_2)$, which is property (b) in Theorem~\ref{ODEthm}. 

We now fix an interval $J\supset[0,\rho_{\max}]$. Before proving that $t_J=\infty$, we will show that 
\begin{equation}\label{varphiJintegralpropertyEq}
\int\varphi_J(t,\rho,x)\,\lambda(dx)=1\qquad\text{for $\rho\in J$ and $t<t_J$.}
\end{equation}
This  then will establishes property (c)  in the statement of Theorem~\ref{ODEthm} for $t\in[0,t_J)$. Then we will use \eqref{varphiJintegralpropertyEq} to derive some estimates on $\varphi_J$ that will yield $\varphi_J>0$ and $t_J=\infty$.

To prove \eqref{varphiJintegralpropertyEq}, we let
$I(t,\rho):=\int\varphi_J(t,\rho,x)\,\lambda(dx)$
 and 
$
\psi_J(t,\rho):=\int x \varphi_J(t,\rho,x)\,\lambda(dx)$.
 We have  
\begin{equation}\label{IODEEq}
 I'(t,\rho)+ \rho I(t,\rho)
  +\psi_J(t,\rho)=\frac{1}{2\bbar\rho}\big(\rho+\psi_J(t,\rho)\big)\Big(\bbar\rho+\int
  x I(t,x)\,\lambda(dx)\Big).
\end{equation}
  This is a (non-homogeneous) affine ODE of the form $I'(t)=b(t)+A(t)I(t)$, where the operator 
$$(A(t)f)(\rho)=-\rho f(\rho)+\frac{1}{2\bbar\rho}\big(\rho+\psi_J(t,\rho)\big)\int
  x f(x)\,\lambda(dx)
$$
is a continuous map from $[0,\delta]$ into the space of bounded linear operators on $C(J)$ for each $\delta<t_J$. Hence this ODE admits a unique solution in $C^1([0,\delta];C(J))$ with initial condition $I(0,\rho)=1$. But~\eqref{IODEEq} is solved by $I(t,\rho)=1$, which  which establishes \eqref{varphiJintegralpropertyEq}.

For the next step, we  let 
$$t_0:=\inf\Big\{t\in[0,t_J)\,\big|\,\min_{\rho_1,\rho_2\in J}\varphi_J(t,\rho_1,\rho_2)<0\Big\}.
$$
Since $\varphi_J$ is a continuous map from $[0,t_J)$ into $C(J\times J)$ and $\varphi_J(0)=1$, we must have $t_0>0$. Due to \eqref{varphiJintegralpropertyEq} we have on  $[0,t_0)$ that
\begin{equation}\label{PreliminaryBound1}
\frac{\rho_1\rho_2}{2\bbar\rho}\le \varphi_J'(t,\rho_1,\rho_2)+(\rho_1+\rho_2)\varphi_J(t,\rho_1,\rho_2)\le\frac{(\rho_1+\rho_{\max})(\rho_2+\rho_{\max})}{2\bbar\rho}.
\end{equation}
When defining 
\begin{equation}\label{hatvarphidef}
\hat\varphi_J(t,\rho_1,\rho_2):=e^{t(\rho_1+\rho_2)}\varphi_J(t,\rho_1,\rho_2),
\end{equation}
the preceding inequality can be rewritten as 
$$\frac{\rho_1\rho_2}{2\bbar\rho}\cdot e^{t(\rho_1+\rho_2)}\le \hat\varphi_J'(t,\rho_1,\rho_2)\le\frac{(\rho_1+\rho_{\max})(\rho_2+\rho_{\max})}{2\bbar\rho}\cdot e^{t(\rho_1+\rho_2)}.
$$
Integrating these inequalities yields that for $0\le t<t_0$ 
   \begin{equation}\label{PreliminaryBound2}
\varphi_J(t,\rho_1,\rho_2)\ge  e^{-t(\rho_1+\rho_2)}+\frac{\rho_1
      \rho_2(1-e^{-t(\rho_1+\rho_2)})}{2\bbar\rho (\rho_1+\rho_2)} >0
\end{equation}
with the convention
$\frac{1-e^{-(\rho_1+\rho_2)t}}{\rho_1+\rho_2}=t$ for $\rho_1=\rho_2=0$. Hence
\begin{equation}\label{PreliminaryBound3}\varphi_J(t,\rho_1,\rho_2)  \le
e^{-t(\rho_1+\rho_2)}+\frac{(\rho_1+\rho_{\max})(\rho_2+\rho_{\max})}{{2\bbar\rho
    (\rho_1+\rho_2)}}
(1-e^{-t(\rho_1+\rho_2)}).
\end{equation}
Inequality \eqref{PreliminaryBound2} ensures that $t_0\ge t_J$. Both inequalities \eqref{PreliminaryBound2} and \eqref{PreliminaryBound3}   ensure the solution $\varphi_J(t)$ does not explode in finite time, which by standard arguments yields  
that $t_J=+\infty$. This proves the global existence of solutions as well as  property (a)  in the statement of Theorem~\ref{ODEthm}. \qed

\bigskip

The preceding lemma works only for measures $\lambda$ that are concentrated on a finite interval. To obtain solutions for more general measures $\lambda$, we need to find upper bounds that are independent of $\rho_{\max}$. 
To this end, we first derive such bounds for the function $\psi(t,\rho)$ defined in \eqref{psiDefinitionEq}. 
By Lemma \ref{RiccatiLemma1}, this function is well-defined whenever $\lambda$ has compact support, and it follows from  dominated convergence together with \eqref{PreliminaryBound2} and \eqref{PreliminaryBound3}  that $\psi\in C^1([0,\infty);C(\bR_+))$ and that 
$\psi'(t,\rho)=\int x \varphi'(t,\rho,x)\,\lambda(dx)$.

\bigskip

\begin{lemma}Under the assumptions of Lemma~\ref{RiccatiLemma1}, we have
\begin{equation}\label{Psibetterbound}
0< \psi(t,\rho)\le \frac{\bbar{\rho^2}}{\bbar{\rho}}\qquad \text{for all $\rho\ge0$.}
\end{equation}
\end{lemma}

\Proof The lower bound in \eqref{Psibetterbound} is clear from $\varphi>0$. To prove the upper bound, 
we suppose by way of contradiction that there exist $t$, $\rho$, and $\eps>0$ such that $\psi(t,\rho)\ge\eps+{\bbar{\rho^2}}/{\bbar{\rho}}$. Then there must be a compact interval $J\supset[0,\rho_{\max}]$ such that
$$\tau_\eps:=\inf\Big\{t\ge0\,\Big|\,\max_{\rho\in J}\psi(t,\rho)\ge\frac{\bbar{\rho^2}}{\bbar{\rho}}+\eps\Big\}
$$
is finite. 
Since $\psi(0,\rho)=\bbar\rho$ and $\bbar\rho^2\le\bbar{\rho^2}$, the time $\tau_\eps$ must also be strictly positive. Moreover, there exists $\rho_\eps\in J$ such that 
$$\max_{\rho\in J}\psi(\tau_\eps,\rho)= \psi(\tau_\eps,\rho_\eps)=\frac{\bbar{\rho^2}}{\bbar{\rho}}+\eps.
$$
Then $\tau_\eps$ is the first time at which the function $t\mapsto \psi(t,\rho_\eps)$ reaches a new maximum, and so
$\psi'(\tau_\eps,\rho_\eps)\ge0
$.

Integrating \eqref{InfiniteDimRiccatiEqn} with respect to $\rho_1\,\lambda(d\rho_1)$ and evaluating at $\rho_2=\rho_\eps$ gives
\begin{eqnarray}\psi'(\tau_\eps,\rho_\eps)+\rho_\eps\psi(\tau_\eps,\rho_\eps)+\int\rho^2\varphi(\tau_\eps,\rho,\rho_\eps)\,\lambda(d\rho)=\frac1{2\bbar\rho}\big(\rho_\eps+\psi(\tau_\eps,\rho_\eps)\big)\Big(\bbar{\rho^2}+\int\rho\psi(\tau_\eps,\rho)\,\lambda(d\rho)\Big).\label{IntegratedRiccatiEq}
\end{eqnarray}
Since $\int\varphi(\tau_\eps,\rho,\rho_\eps)\,\lambda(d\rho)=1$, the Cauchy--Schwarz inequality (or, alternatively, Jensen's inequality) implies that
$\int\rho^2\varphi(\tau_\eps,\rho,\rho_\eps)\,\lambda(d\rho)\ge\psi(\tau_\eps,\rho_\eps)^2$.
Moreover, the definition of $\rho_\eps$ and the fact that $\lambda$ is supported on  $J$  yield that
$\int\rho\psi(\tau_\eps,\rho)\,\lambda(d\rho)\le\bbar\rho\psi(\tau_\eps,\rho_\eps)$.
Plugging these two inequalities into \eqref{IntegratedRiccatiEq} leads to
\begin{eqnarray*}\psi'(\tau_\eps,\rho_\eps)&\le&\frac1{2\bbar\rho}\big(\rho_\eps+\psi(\tau_\eps,\rho_\eps)\big)\big(\bbar{\rho^2}+\bbar\rho\psi(\tau_\eps,\rho_\eps)\big)-\rho_\eps\psi(\tau_\eps,\rho_\eps)-\psi(\tau_\eps,\rho_\eps)^2\\
&=&\frac{\rho_\eps\bbar{\rho^2}}{2\bbar\rho}+\Big(\frac{\bbar{\rho^2}}{2\bbar\rho}-\frac{\rho_\eps}2\Big)\psi(\tau_\eps,\rho_\eps)-\frac12\psi(\tau_\eps,\rho_\eps)^2=: p(\psi(\tau_\eps,\rho_\eps)),
\end{eqnarray*}
where $p(\cdot)$ is a polynomial function of degree two. 
It has the two roots $-\rho_\eps\le0$ and 
${\bbar{\rho^2}}/{\bbar\rho}>0$.
Therefore $p(x)<0$ for $x>{\bbar{\rho^2}}/{\bbar\rho}$ and in turn $0>p(\psi(\tau_\eps,\rho_\eps))=\psi'(\tau_\eps,\rho_\eps)$,
which contradicts the fact that $\psi'(\tau_\eps,\rho_\eps)\ge0
$. \qed

\bigskip

\bigskip

\begin{lemma}\label{betterboundLemma}
Under the assumptions of Lemma~\ref{RiccatiLemma1}, we have
\begin{equation}\label{betterboundEq}
\begin{split} e^{-t(\rho_1+\rho_2)}+&\frac{\rho_1
      \rho_2(1-e^{-t(\rho_1+\rho_2)})}{2\bbar\rho (\rho_1+\rho_2)}\le\varphi(t,\rho_1,\rho_2)\\ & \le
\exp(-(\rho_1+\rho_2)t)+\frac{(\rho_1+\frac{\bbar{\rho^2}
  }{\bbar{\rho}})(\rho_2+\frac{\bbar{\rho^2} }{\bbar{\rho}} )}{{2\bbar\rho (\rho_1+\rho_2)}}
\big(1-\exp(-(\rho_1+\rho_2)t)\big),
\end{split}
\end{equation}
\begin{equation}\label{phi'betterboundEq}
-(\rho_1+\rho_2)-\frac{(\rho_1+\frac{\bbar{\rho^2} }{\bbar{\rho}})(\rho_2+\frac{\bbar{\rho^2} }{\bbar{\rho}})}
    {{2\bbar\rho}} \le\varphi'(t,\rho_1,\rho_2)\le \frac{(\rho_1+\frac{\bbar{\rho^2} }{\bbar{\rho}})(\rho_2+\frac{\bbar{\rho^2} }{\bbar{\rho}})}
    {{2\bbar\rho}}.
\end{equation}
\end{lemma}

\Proof The ODE~\eqref{InfiniteDimRiccatiEqn} can be rewritten as
\begin{equation}\label{AltODEEq}
\varphi'(t,\rho_1,\rho_2)+(\rho_1+\rho_2)\varphi(t,\rho_1,\rho_2)=\frac1{2\bbar\rho}\big(\rho_1+\psi(t,\rho_1)\big)\big(\rho_2+\psi(t,\rho_2)\big). 
\end{equation}
Defining $\hat\varphi$ as in 
\eqref{hatvarphidef} and using the upper bound in~\eqref{Psibetterbound} thus yields that
\begin{equation}\label{betterVarPhiHatBoundEq}
\frac{\rho_1\rho_2}{2\bbar\rho}\cdot e^{t(\rho_1+\rho_2)}\le \hat\varphi'(t,\rho_1,\rho_2)\le \frac{(\rho_1+\frac{\bbar{\rho^2} }{\bbar{\rho}})(\rho_2+\frac{\bbar{\rho^2} }{\bbar{\rho}})}
    {{2\bbar\rho}}\cdot e^{t(\rho_1+\rho_2)}.
\end{equation}
Arguing as in the final step of the proof Lemma~\ref{RiccatiLemma1} now yields~\eqref{betterboundEq}. By plugging~\eqref{betterboundEq} back into~\eqref{AltODEEq} and using once again \eqref{Psibetterbound}, we obtain~\eqref{phi'betterboundEq}.\qed

\bigskip

\begin{lemma}\label{LipschitzLemma}For all $R$, $T>0$ there exist constants $L_1$, $L_2\ge0$ depending only on $R$, $T$, $\bbar\rho$, and $\bbar{\rho^2}$ such that for all $t\in[0,T]$ and $\rho_1,\wt\rho_1\,\rho_2\in[0,R]$,
\begin{equation}\label{phiLipEq}
|\varphi(t,\rho_1,\rho_2)-\varphi(t,\wt\rho_1,\rho_2)|\le L_1|\rho_1-\wt\rho_1|\quad\text{}
\end{equation}
and
\begin{equation}\label{phiprimeLipEq}
|\varphi'(t,\rho_1,\rho_2)-\varphi'(t,\wt\rho_1,\rho_2)|\le L_2|\rho_1-\wt\rho_1|.
\end{equation}
\end{lemma}

\Proof We consider $\rho_1,\wt\rho_1,\rho_2\ge 0$ and define 
$$\Delta\rho_1:=\wt\rho_1-\rho_1,\quad \Delta\varphi(t)=\varphi(t,\wt\rho_1,\rho_2)-\varphi(t,\rho_1,\rho_2)\quad\text{and}\quad \Delta\psi(t)=\psi(t,\wt\rho_1)-\psi(t,\rho_1).$$
 By subtracting the equation \eqref{AltODEEq}
satisfied by $\varphi(t,\rho_1,\rho_2)$ from the corresponding one satisfied by
$\varphi(t,\wt\rho_1,\rho_2)$, we get
\begin{equation}\label{psiPrimeLipIneq}
\Delta\varphi'(t)+\varphi(t,\wt\rho_1,\rho_2)\Delta\rho_1+(\rho_2+\rho_1)\Delta\varphi(t)=\frac1{2\bbar\rho}\Big(\rho_2+\psi(t,\rho_2)\Big)\Big(\Delta\rho_1+\Delta\psi(t)\Big).
\end{equation}
 This equation is a linear non-homogeneous ODE for $\Delta\varphi(t)$ and, since $\Delta\varphi(0)=0$, solved by
\begin{eqnarray*}\lefteqn{\Delta\varphi(t)=}\\
&&\int_0^t\left[\Big(\frac1{2\bbar\rho}
(\rho_2+\psi(s,\rho_2))-\varphi(s,\wt\rho_1,\rho_2)\Big)\Delta\rho_1+\frac1{2\bbar\rho}\Big(\rho_2+\psi(s,\rho_2)\Big)\Delta\psi(s)\right]e^{-(\rho_1+\rho_2)(t-s)}\,ds,\end{eqnarray*}
 Since $|\psi(s,\rho_2)|\le
{\bbar{\rho^2}}/{\bbar{\rho}}$, we get with \eqref{betterboundEq} and
$\sup_{\alpha\ge 0} \frac{1-e^{-\alpha t}}{\alpha}=t$ that
\begin{equation}\label{ub_deltaphi}|\Delta\varphi(t)|\le
\frac1{2\bbar\rho}\Big(\rho_2+\frac{\bbar{\rho^2}}{\bbar{\rho}}\Big) \int_0^t\big(|\Delta\rho_1|+|\Delta\psi(s)|\big)\,ds+\bigg(1+\frac{(\wt\rho_1+\frac{\bbar{\rho^2} }{\bbar{\rho}})(\rho_2+\frac{\bbar{\rho^2} }{\bbar{\rho}} )}{{2\bbar\rho}}T\bigg)\int_0^t|\Delta\rho_1|\,ds.
\end{equation}

Now, we have that
\begin{eqnarray*}|\Delta \psi(t)|&=&\Big|\int
(\varphi(t,\wt\rho_1,x)-\varphi(t,\rho_1,x))x\,\lambda(dx) \Big|  
\le\int |\varphi(t,\wt\rho_1,x)-\varphi(t,\rho_1,x)|x\,\lambda(dx)  \\
&\le&\frac{\bbar{\rho^2}}{\bbar\rho}\int_0^t\big(|\Delta\rho_1|
+|\Delta\psi(s)|\big)\,ds+\bigg(\bbar\rho+\Big(\wt\rho_1+\frac{\bbar{\rho^2} }{\bbar{\rho}}\Big)\frac{\bbar{\rho^2} }{\bbar{\rho}}T\bigg)\int_0^t|\Delta\rho_1|\,ds.
\end{eqnarray*}
For the last inequality, we have used Fubini's theorem
and~\eqref{ub_deltaphi}. Now, Gronwall's Lemma gives:
\begin{equation}\label{DeltaPsiIneq}
  |\Delta \psi(t)|\le |\Delta\rho_1| t\bigg(\frac{\bbar{\rho^2}}{\bbar\rho} +\bbar\rho+ \Big(\wt\rho_1+\frac{\bbar{\rho^2} }{\bbar{\rho}}\Big)\frac{\bbar{\rho^2} }{\bbar{\rho}}T\bigg)\exp\Big(\frac{\bbar{\rho^2}}{\bbar\rho}t\Big).
\end{equation}
Plugging this back into~\eqref{ub_deltaphi}, we get the existence of a constant $L_1$, which depends only on $R$, $T$, $\bbar\rho$, and $\bbar{\rho^2}$, such that 
\begin{equation}\label{DeltaPhiIneq}
|\Delta\varphi(t)|\le
L_1 |\Delta\rho_1|
 . 
\end{equation}
Finally, using \eqref{DeltaPsiIneq} and \eqref{DeltaPhiIneq} in \eqref{psiPrimeLipIneq} and recalling the locally uniform bounds \eqref{betterboundEq} and \eqref{Psibetterbound} on $\varphi$ and $\psi$ gives \eqref{phiprimeLipEq}. \qed

\bigskip

Now drop the assumption that $\lambda$ is supported on a compact interval and aim at proving existence and uniqueness of solutions in this general case. To this end, we  take a sequence $R_n\ua\infty$ for which $\lambda([0,R_1])>0$ and define 
\begin{equation}\label{lambdanApproximationEq}
f_n:=\frac1{\lambda([0,R_n])}\Ind{[0,R_n]}\qquad\text{and}\qquad d\lambda_n=f_n\,d\lambda,
\end{equation}
so that each $\lambda_n$ satisfies the assumptions of Lemma~\ref{RiccatiLemma1}. By $\varphi_n$ we denote the corresponding solution of ~\eqref{InfiniteDimRiccatiEqn},~\eqref{InfiniteDimRiccatiEqnInitialCond}  provided by that lemma.
For each $n\ge1$, we have
\begin{equation}\label{rhonrho2nEq}
\bbar\rho_n:=\int \rho\,\lambda_n(d\rho)\ge \int_0^{R_1}\rho\,\lambda(d\rho) =:\bbar\rho_0
\quad\text{and}\quad \bbar{\rho^2_n}:=\int\rho^2\,\lambda_n(d\rho)\le \frac{\bbar{\rho^2}}{\lambda([0,R_1])}=:\bbar{\rho^2_0}.
\end{equation}
Hence, Lemma~\ref{betterboundLemma} yields that for each $n$, \begin{equation}\label{n betterboundEq}
e^{-t(\rho_1+\rho_2)}+\frac{\rho_1
      \rho_2(1-e^{-t(\rho_1+\rho_2)})}{2\bbar\rho_0 (\rho_1+\rho_2)}\le\varphi_n(t,\rho_1,\rho_2)  \le
1+\frac{(\rho_1+\frac{\bbar{\rho^2_0}
  }{\bbar{\rho}_0})(\rho_2+\frac{\bbar{\rho^2_0} }{\bbar{\rho}_0}
  )}{{2\bbar\rho_0}} \frac{1-e^{-(\rho_1+\rho_2)t}}{\rho_1+\rho_2}
\end{equation}
and 
\begin{equation}\label{n phi'betterboundEq}
-(\rho_1+\rho_2)-\frac{(\rho_1+\frac{\bbar{\rho^2_0} }{\bbar{\rho}_0})(\rho_2+\frac{\bbar{\rho^2_0} }{\bbar{\rho}_0})}
    {{2\bbar\rho_0}}\le\varphi'_n(t,\rho_1,\rho_2)\le \frac{(\rho_1+\frac{\bbar{\rho^2_0} }{\bbar{\rho}_0})(\rho_2+\frac{\bbar{\rho^2_0} }{\bbar{\rho}_0})}
    {{2\bbar\rho_0}}.
\end{equation}
Similarly, Lemma \ref{LipschitzLemma} yields that for all $R$, $T>0$ there is a constant  $L\ge0$ such that for all $n$
\begin{equation}\label{phiLipEq2}
|\varphi_n(t,\rho_1,\rho_2)-\varphi_n(t,\wt\rho_1,\rho_2)|\le L|\rho_1-\wt\rho_1|\quad\text{for all $t\in[0,T]$, $\rho_1,\wt\rho_1,\rho_2\in[0,R]$,}
\end{equation}

The inequalities \eqref{n phi'betterboundEq}, \eqref{phiLipEq2} and the Arzela--Ascoli theorem imply that the sequence $(\varphi_n)_{n\in\bN}$ is relatively compact in the class of continuous functions on $[0,T]\times[0,R]^2$ whenever $T$, $R>0$, and hence admits a convergent subsequence in that class. By passing to a subsequence arising from a diagonalization argument  if necessary, we may assume  that there exists a continuous function $\varphi:[0,\infty)\times\bR_+^2\to\bR_+$ such that $\varphi_{n}\to\varphi$ locally uniformly. 

The uniform bounds \eqref{n betterboundEq} and dominated convergence imply that 
\begin{eqnarray}\label{psiNconvergenceEq}
\psi_n(t,\rho)&:=&\int x\varphi_n(t,\rho,x)\,\lambda_n(dx)=\int x\varphi_n(t,\rho,x)f_n(x)\,\lambda(dx)\\
&&\longrightarrow \int x\varphi(t,\rho,x)\,\lambda(dx)=\psi(t,\rho)
\end{eqnarray}
locally uniformly in $(t,\rho)$. 
 Hence, 
$F_{\lambda_n}(\varphi_n(t))(\rho_1,\rho_2)\to F_{\lambda}(\varphi(t))(\rho_1,\rho_2)
$, locally uniformly in $(t,\rho_1,\rho_2)$, where $F_{\lambda_n}$ is defined through \eqref{F(f)Eq}.  
Since $\varphi_n'=F_{\lambda_n}(\varphi_n)$, we conclude that $\varphi_n'\to F_\lambda(\varphi)$ locally uniformly in $[0,\infty)\times\bR^2_+$. Moreover, we have for each $n$ that
$$\varphi_n(t,\rho_1,\rho_2)-1=\int_0^t\varphi'_n(s,\rho_1,\rho_2)\,ds.
$$
The left-hand side of this equation converges to $\varphi(t,\rho_1,\rho_2)-1$, whereas the right-hand side converges to $\int_0^tF_\lambda(\varphi(s))(\rho_1,\rho_2)\,ds$. This proves that $\varphi$ solves \eqref{InfiniteDimRiccatiEqn} and that $\varphi\in C^1([0,\infty);C(\bR^2_+))$.

\begin{remark}By  sending $R_1$ to infinity in \eqref{n betterboundEq} we get that the solution $\varphi$ constructed above satisfies the bounds
\begin{equation}\label{varphiFinalEstimates}
 e^{-t(\rho_1+\rho_2)}+\frac{\rho_1
      \rho_2(1-e^{-t(\rho_1+\rho_2)})}{2\bbar\rho(\rho_1+\rho_2)}\le\varphi(t,\rho_1,\rho_2)\le 1+\frac{(\rho_1+\frac{\bbar{\rho^2} }{\bbar{\rho}})(\rho_2+\frac{\bbar{\rho^2} }{\bbar{\rho}} )}{{2\bbar\rho}}\frac{1-e^{-(\rho_1+\rho_2)t}}{\rho_1+\rho_2}.
\end{equation}
 From \eqref{Psibetterbound}, \eqref{psiNconvergenceEq}, and the lower bound in  \eqref{varphiFinalEstimates} we get moreover that 
\begin{equation}\label{psiFinalEstimate}
0<\psi(t,\rho)\le\frac{\bbar{\rho^2}}{\bbar\rho}.
\end{equation}
\end{remark}

\bigskip

Now we turn to prove the uniqueness of solutions in the class of functions $\varphi\in C^1([0,\infty);C(\bR^2_+))$ satisfying a bound of the form \eqref{phiiEstimateEq}
To this end, let $\varphi_1$ and $\varphi_2$ be two solutions in that class and set 
$$\delta(t,\rho_1,\rho_2)=\varphi_2(t,\rho_1,\rho_2)-\varphi_1(t,\rho_1,\rho_2).$$ 
We will show that $\| \delta(t) \|_{L^2(\wt\lambda
\otimes\wt\lambda)}^2=0$ for all $t$, whenever  $\wt\lambda$ is a positive finite Borel measure  of the form
$\wt\lambda=\lambda+\mu$, 
where $\mu$ is a positive finite Borel measure with compact support. Taking, for instance, $\mu$  as the Lebesgue measure on $[0,R]$ will then imply that  $\delta(t,\rho_1,\rho_2)=0$ for $\rho_1,\rho_2\in[0,R]$. So this will give the uniqueness of solutions.  

Let us define $\wt F_\lambda$ as
\begin{equation}\label{tildeF(f)Eq}
\wt F_\lambda(f)(\rho_1,\rho_2)=\frac1{2\bbar\rho}\Big(\rho_1+\int x f(\rho_1,x)\,\lambda(dx)\Big)\Big(\rho_2+\int xf(x,\rho_2)\,\lambda(dx)\Big).
\end{equation}

\bigskip

\begin{lemma}\label{tildeF_lip_lemma}We have $\wt F_\lambda(\varphi_i(t))\in L^2(\wt\lambda\otimes\wt\lambda)$ and 
\begin{equation}\label{tildeF_lip}\|\wt F_\lambda(\varphi_1(t))-\wt F_\lambda(\varphi_2(t))\|_{L^2(\wt\lambda
\otimes\wt\lambda)}^2 \le C(\|\delta(t)\|_{L^2(\wt\lambda
\otimes\wt\lambda)}^2+\|\delta(t)\|_{L^2(\wt\lambda
\otimes\wt\lambda)}^4 ),
\end{equation}
where  $C$ is a positive constant that depends only on
$\bbar{\rho}$ and $\bbar{\rho^2}$.
\end{lemma}

\Proof For simplicity, we will drop the argument $t$ throughout the proof. We may write
\begin{eqnarray*}
\lefteqn{\wt F_\lambda(\varphi_2)-\wt F_\lambda(\varphi_1)}\\&=& \frac1{2\bbar\rho}\bigg[\int x
\delta(\rho_1,x)\,\lambda(dx)\Big(\rho_2+\int x\varphi_1(x,\rho_2)\,\lambda(dx)\Big)\\
&&
+\int x\delta(x,\rho_2)\,\lambda(dx)\Big(\rho_1+\int x
\varphi_1(\rho_1,x)\,\lambda(dx)\Big) + \int x
\delta(\rho_1,x)\,\lambda(dx)\int x\delta(x,\rho_2)\,\lambda(dx)\bigg].
\end{eqnarray*}
Thus,
\begin{eqnarray*}
\lefteqn{(\wt F_\lambda(\varphi_2)-\wt F_\lambda(\varphi_1))^2}\\
&\le&\frac3{4\bbar\rho^2}\Bigg[ \left(\int x
\delta(\rho_1,x)\,\lambda(dx)\Big(\rho_2+\int x\varphi_1(x,\rho_2)\,\lambda(dx)\Big)\right)^2\\
&&
+\left(\int x\delta(x,\rho_2)\,\lambda(dx)\Big(\rho_1+\int x
\varphi_1(\rho_1,x)\,\lambda(dx)\Big)\right)^2 +  \left(\int x
\delta(\rho_1,x)\,\lambda(dx)\int x\delta(x,\rho_2)\,\lambda(dx)\right)^2\Bigg].
\end{eqnarray*}
Now we integrate this inequality with respect to $\wt\lambda(d\rho_1)\,\wt
\lambda(d\rho_2)$. The two first terms can be analyzed in the same way. First, we
observe that $\int
(\rho_2+\frac{\bbar{\rho^2}}{\bbar{\rho}})^2\,\wt\lambda(d\rho_2) $ is finite.  
Then we note that 
\begin{equation}\label{intermed_uniq}\left( \int x\delta(\rho_1,x)\, \lambda(dx) \right)^2\le \bbar{\rho^2}  \int
\delta(\rho_1,x)^2 \,\lambda(dx)\le  \bbar{\rho^2}  \int
\delta(\rho_1,x)^2 \,\wt\lambda(dx).
\end{equation}
Hence,
$$\int\left( \int x\delta(\rho_1,x)\, \lambda(dx) \right)^2\,\wt\lambda(d\rho_1)\le\bbar{\rho^2}\|\delta\|_{L^2(\wt\lambda
\otimes\wt\lambda)}^2.
$$
Thus, the two first terms can be bounded by $C_0\|\delta\|_{L^2(\wt\lambda
\otimes\wt\lambda)}^2$, where $C_0$ is a constant that only depends on
$\bbar{\rho}$ and $\bbar{\rho^2}$. Using once again~\eqref{intermed_uniq}, we
get that the third term can be  bounded from above by $C_1\|\delta\|_{L^2(\wt\lambda
\otimes\wt\lambda)}^4$, where the constant $C_1$ depends only on
$\bbar{\rho}$ and $\bbar{\rho^2}$.\qed

\bigskip

Now we differentiate $\delta^2$ and integrate over $[0,t]$:
\begin{eqnarray*} \delta(t,\rho_1,\rho_2)^2 &=& -2\int_0^t(\rho_1+\rho_2)
\delta(s,\rho_1,\rho_2)^2\, ds\\&& +
2\int_0^t\delta(s,\rho_1,\rho_2)\Big[\wt F_\lambda(\varphi_2(s ))(\rho_1,\rho_2)-
\wt F_\lambda(\varphi_1(s ))(\rho_1,\rho_2)\Big]\,ds \\
&\le&2\int_0^t\delta(s,\rho_1,\rho_2)\Big[\wt F_\lambda(\varphi_2(s ))(\rho_1,\rho_2)-
\wt F_\lambda(\varphi_1(s ))(\rho_1,\rho_2)\Big]\,ds.
\end{eqnarray*}
We now integrate w.r.t.  $\wt\lambda(d\rho_1)\,
\wt\lambda(d\rho_2)$ and get by using the Cauchy--Schwarz inequality,
$$\| \delta(t,\cdot ) \|_{L^2(\wt\lambda
\otimes\wt\lambda)}^2 \le 2 \int_0^t \| \delta(s,\cdot ) \|_{L^2(\wt\lambda
\otimes\wt\lambda)} \| \wt F_\lambda(\varphi_2(s ))-\wt F_\lambda(\varphi_1(s )) \|_{L^2(\wt\lambda
\otimes\wt\lambda)}\,ds.$$
By continuity of~$t\rightarrow \| \delta(t,\cdot) \|_{L^2(\wt\lambda
\otimes\wt\lambda)}$, we know that for each $T>0$ there is a constant $K$ such that $\| \delta(t,\cdot) \|_{L^2(\wt\lambda
\otimes\wt\lambda)}\le K$ when $t\in[0,T]$. Thus,  we get from Lemma \ref{tildeF_lip_lemma}
that
$$\| \delta(t,\cdot) \|_{L^2(\wt\lambda
\otimes\wt\lambda)}^2 \le \sqrt{C (1+K^2)}\int_0^t \| \delta(s,\cdot) \|_{L^2(\wt\lambda
\otimes\wt\lambda)} ^2\,ds,$$
which in turn gives that $\| \delta(t,\cdot) \|_{L^2(\wt\lambda
\otimes\wt\lambda)}^2=0$ on $[0,T]$ by  Gronwall's Lemma. This concludes the proof of uniqueness.

\bigskip

Now we turn to proving the properties (a) through (f) in Theorem \ref{ODEthm}. Property (a) (strict positivity) can be proved just as in the case of a compactly supported measure $\lambda$ in Lemma \ref{RiccatiLemma1}. Property (b) (symmetry) is already clear. Property {(c)} ($\int\varphi(t,\rho,x)\,\lambda(dx)=1$) follows from the corresponding property of the approximating functions $\varphi_n$, the uniform bounds \eqref{n betterboundEq}, and dominated convergence. 

Property (d) states that   $\varphi\in C^{2}([0,\infty);C(\bR^2_+))$. By
dominated convergence and the bound \eqref{n phi'betterboundEq}, which also
holds for $\varphi'$ in place of $\varphi_n'$, we get that $\psi(t,\rho)$ belongs to $ C^{1}([0,\infty);C(\bR_+))$. Thus, our ODE gives $\varphi'\in  C^{1}([0,\infty);C(\bR^2_+))$, which proves property (d).

We now prove property  (e). It is clearly enough to prove it when $f:\bR_+\to\bR$ is a bounded measurable function with compact support. To this end, let
$\hat\varphi(t,\rho_1,\rho_2):=e^{t(\rho_1+\rho_2)}\varphi(t,\rho_1,\rho_2)$.
Then $\hat\varphi$ belongs to $C^1([0,\infty);C(\bR^2_+))$ and 
$$\hat\varphi'(t,\rho_1,\rho_2)=\frac1{2\bbar\rho}\cdot e^{t(\rho_1+\rho_2)}(\rho_1+\psi(t,\rho_1))(\rho_2+\psi(t,\rho_2)).
$$
That is, $\hat\varphi'(t,\rho_1,\rho_2)=g(t,\rho_1)g(t,\rho_2)$ for a function $g$. Thus, 
$$\int\int f(x_1)f(x_2)\hat\varphi'(t,x_1,x_2)\,\lambda(dx_1)\,\lambda(dx_2)=\Big(\int f(x)g(t,x)\,\lambda(dx)\Big)^2\ge0.
$$
Since
$\hat\varphi(t)=1+\int_0^t\hat\varphi'(s)\,ds$,
we find that $\hat\varphi(t)$ is nonnegative definite. Finally, with $\hat f(x)=e^{-t x}f(x)$,
$$\int\int f(x_1)f(x_2)\varphi(t,x_1,x_2)\,\lambda(dx_1)\,\lambda(dx_2)=\int\int \hat f(x_1)\hat f(x_2)\hat\varphi(t,x_1,x_2)\,\lambda(dx_1)\,\lambda(dx_2)\ge0.
$$
This establishes property  (e) in Theorem \ref{ODEthm}.

Finally, property (f) (the local Lipschitz property for $\varphi$ and $\varphi'$) follows just as in Lemma \ref{LipschitzLemma}. This concludes the proof of Theorem \ref{ODEthm}.
\qed

\subsection{Proof of Theorem \ref{verificationThm}}

The strategy in the proof of Theorem \ref{verificationThm} is to use a verification argument and based on guessing the optimal costs $V(T,E(\cdot),x)$ for liquidating $x$ shares over $[0,T]$ with additional and arbitrary initial data $E(\cdot)$. The result of our guess is formula 
\eqref{VdefEq} below. We explain its heuristic derivation in Appendix~\ref{app_heur}. 

Let $\varphi$ be a solution of the infinite-dimensional Riccati equation ~\eqref{InfiniteDimRiccatiEqn},~\eqref{InfiniteDimRiccatiEqnInitialCond}. This solution gives rise to  a family of linear operators $\Phi_t:L^2(\lambda)\to L^2(\lambda)\cap C(\bR_+)$ defined by 
$$\Phi_tf(\rho)=\int f(x)\varphi(t,x,\rho)\,\lambda(dx),\qquad f\in L^2(\lambda).
$$
By  \eqref{phiiEstimateEq}, $t\mapsto \Phi_tf$ is a continuous map into both $L^2(\lambda)$ and $C(\bR_+)$ for each $f\in L^2(\lambda)$. By the inequality \eqref{n phi'betterboundEq}, which also hold for $\varphi$ in place of $\varphi_n$, $t\mapsto \Phi_tf$ is a continuously differentiable map into both $L^2(\lambda)$ and $C(\bR_+)$ for each $f\in L^\infty(\lambda)$.

For $t\ge0$, $E(\cdot)\in L^2(\lambda)$, and $x\in\bR$, we define 
\begin{eqnarray}\label{VdefEq}
{V}(t,E(\cdot),x):=\frac12\Big[\frac1{\varphi_0(t)}\big(x-(\Phi_tE)(0)\big)^2-\langle E,\Phi_tE\rangle\Big],
\end{eqnarray}
where $\<\cdot,\cdot\>$ denotes the usual inner product in $L^2(\lambda)$. For $t\in[0,T]$ and a $[0,T]$-admissible strategy $X$ we define
$$C^X_t:=\int_{[0,t)}\int E^X_s\,d\lambda\,dX_s+\frac12\sum_{s<t}(\Delta X_s)^2+{V}(T-t,E_t^X,X_t).
$$
By Lemma \ref{CostFunctionalLemma}, the first two terms on the right correspond to the cost accumulated by the strategy up to time $t$. Moreover,
$${V}(0,E_T^X,X_T)=\frac12(X_T)^2-\int E^X_T\,d\lambda\cdot X_T=\frac12(\Delta X_T)^2+\int E_T^X\,d\lambda\cdot\Delta X_T,
$$ 
due to the requirement $X_{T+}=0$. This gives 
$C_T^X=\cC(X)$.
Our goal is thus to show the following verification lemma:  $dC^X_t\ge0$ with equality if and only if $X=X^*$ for a certain strategy $X^*$. This will identify $X^*$ as the optimal strategy and ${V}(T,E(\cdot),x)$ as the optimal cost for liquidating $x$ shares over $[0,T]$ with additional initial data $E(\cdot)$ at time $t=0$. In the formalism of potential theory, ${V}(T,0,-1)$ will then be  the minimal energy of a probability measure on $[0,T]$.

\bigskip

\begin{lemma}For every $[0,T]$-admissible strategy $X$, $C^X_t$ is absolutely continuous in $t$ and 
\begin{equation}\label{dCXdtEq}
\frac{dC^X_t}{dt}=\frac12\bigg[\frac{\psi(T-t,0)}{\varphi_0(T-t)}\big(X_t-(\Phi_{T-t}E^X_t)(0)\big)+\int E^X_t(\rho)\big(\rho+\psi(T-t,\rho)\big)\,\lambda(d\rho)\bigg]^2
\end{equation}
for a.e. $t\in[0,T]$.
\end{lemma}

\Proof Recall the following integration by parts formula for left-continuous functions $\alpha_t,\,\beta_t$ of locally bounded variation:
$$\alpha_t\beta_t-\alpha_s\beta_s=\int_{[s,t)}\alpha_r\,d\beta_r+\int_{[s,t)}\beta_r\,d\alpha_r+\sum_{r\in[s,t)}\Delta\alpha_r\Delta\beta_r.
$$
It follows that $t\mapsto E^X_t(\rho)$ is of bounded variation and
\begin{equation}\label{dEtEquation}
E^X_t(\rho)-E^X_s(\rho)=X_t-X_s-\rho \int_s^tE^X_r(\rho)\,dr
\end{equation}
as well as
\begin{equation}\label{dEtEquation2}\begin{split}
\lefteqn{E^X_t(\rho_1)E^X_t(\rho_2)-E^X_s(\rho_1)E^X_s(\rho_2)}\\
&=\int_{[s,t)}\big(E^X_r(\rho_1)+E^X_r(\rho_2)\big)\,dX_r-\int_s^t(\rho_1+\rho_2)E^X_r(\rho_1)E^X_r(\rho_2)\,dr+\sum_{r\in[s,t)}(\Delta X_r)^2.
\end{split}
\end{equation}
Therefore,
\begin{eqnarray*}\lefteqn{\varphi(T-t,\rho_1,\rho_2)E^X_t(\rho_1)E_t^X(\rho_2)-\varphi(T-s,\rho_1,\rho_2)E^X_s(\rho_1)E_s^X(\rho_2)}\\
&&=-\int_s^t\varphi'(T-r,\rho_1,\rho_2)E^X_r(\rho_1)E_r^X(\rho_2)\,dr-\int_s^t\varphi(T-r,\rho_1,\rho_2)E^X_r(\rho_1)E_r^X(\rho_2)(\rho_1+\rho_2)\,dr\\
&&\qquad+\int_{[s,t)}\varphi(T-r,\rho_1,\rho_2)\big(E^X_r(\rho_1)+E_r^X(\rho_2)\big)\,dX_r+\sum_{r\in[s,t)}\varphi(T-r,\rho_1,\rho_2)(\Delta X_r)^2.
\end{eqnarray*}
We have already observed in the proof of Lemma \ref{CostFunctionalLemma} that  $|E_r^X(\rho)|$ is uniformly bounded in $r\in[0,T]$ and $\rho\ge0$ by the total variation of $X$.   Hence we may integrate both sides of the preceding identity with respect to $\lambda(d\rho_1)\,\lambda(d\rho_2)$ to obtain, with the symmetry of $\varphi$ and the notation $\wh E^X_t:=\rho E^X_t(\rho)$, that
\begin{eqnarray*}\lefteqn{\<E^X_t,\Phi_{T-t} E^X_t\>-\<E^X_s,\Phi_{T-s} E^X_s\>}\\
&&=-\int_s^t\<E^X_r,\Phi'_{T-r} E_r^X\>\,dr-2\int_s^t\<\wh E^X_r,\Phi_{T-r} E^X_r\>\,dr+2\int_{[s,t)}\<1, E^X_r\>\,dX_r+\sum_{r\in[s,t)}(\Delta X_r)^2.
\end{eqnarray*}
By a similar reasoning  we obtain
$$\Phi_{T-t}E^X_t(0)-\Phi_{T-s}E^X_s(0)=-\int_s^t\Phi'_{T-r}E^X_r(0)\,dr-\int_s^t\Phi_{T-r}\wh E^X_r(0)\,dr+X_t-X_s
$$
and
$$\big(X_t-\Phi_{T-t}E^X_t(0)\big)^2-\big(X_s-\Phi_{T-s}E^X_s(0)\big)^2=2\int_s^t\big(X_r-\Phi_{T-r}E^X_r(0)\big)\big(\Phi'_{T-r}E^X_r(0)+\Phi_{T-r}\wh E^X_r(0)\big)\,dr.
$$
Using these formulas, we can now compute
\begin{eqnarray*}{C^X_t-C^X_s}
&=&\frac12\int_s^t\frac{\varphi_0'(T-r)}{\varphi_0^2(T-r)}\big(X_r-\Phi_{T-r}E^X_r(0)\big)^2\,dr\\
&&\qquad+\int_s^t\frac{X_r-\Phi_{T-r}E^X_r(0)}{\varphi_0(T-r)}\big(\Phi'_{T-r}E^X_r(0)+\Phi_{T-r}\wh E^X_r(0)\big)\,dr\\
&&\qquad+ \frac12\int_s^t\<E^X_r,\Phi'_{T-r} E_r^X\>\,dr+\int_s^t\<\wh E^X_r,\Phi_{T-r} E^X_r\>\,dr
\end{eqnarray*}
Therefore, $C^X_t$ is absolutely continuous on $[0,T]$ and has the derivative
\begin{align*}\frac{dC^X_t}{dt}&=\frac12\frac{\varphi_0'(T-t)}{\varphi_0^2(T-t)}\big(X_t-\Phi_{T-t}E^X_t(0)\big)^2\\
&
\qquad+\frac{X_t-\Phi_{T-t}E^X_t(0)}{\varphi_0(T-t)}\big(\Phi'_{T-t}E^X_t(0)+\Phi_{T-t}\wh E^X_t(0)\big)+\frac12\<E^X_t,\Phi'_{T-t} E_t^X\>+\<\wh E^X_t,\Phi_{T-t} E^X_t\>
\end{align*}
for a.e. $t$.

To further analyze the preceding formula, we take an extra point $\Delta$. We let $\bbar \lambda:=\lambda+\delta_\Delta$ and extend $E^X_t$ and $\varphi$ to functions on $\{\Delta\}\cup[0,\infty)$ by putting
\begin{equation}\label{def_EX_extended}
\begin{split}
E^X_t(\Delta)&:=\frac1{\varphi_0(T-t)}\big(X_t-(\Phi_{T-t} E^X_t)(0)\big),\\
t\ge 0, \varphi(t,\Delta,\rho)=\varphi(t,\rho,\Delta)&:=\varphi(t,0,\rho),\\\varphi(t,\Delta,\Delta)&:=\varphi(t,0,0)=\varphi_0(t).
\end{split}
\end{equation}
We furthermore define the function
$$f(x)=\begin{cases}x&\text{if $x\in[0,\infty)$,}\\
0&\text{if $x=\Delta$,}
\end{cases}
$$
and we extend the definition of $\wh E^X$ via $\wh E^X_t(x)=f(x)E^X_t(x)$
for $x\in\{\Delta\}\cup[0,\infty)$. Finally, we set for $g \in
L^2(\bbar{\lambda})$, $\bbar{\Phi}_t g(x)=\int
\varphi(t,x,y)g(y)\bbar{\lambda}(dy)$ and one easily  checks that
$\bbar{\Phi}:L^2(\bbar{\lambda})\rightarrow L^2(\bbar{\lambda})$.

With this notation, we get
\begin{eqnarray*}
\frac{dC^X_t}{dt}&=&\frac12\<E^X_t,\bbar{\Phi}'_{T-t} E^X_t\>_{L^2(\bbar\lambda)}+\<E^X_t,\bbar{\Phi}_{T-t} \wh E^X_t\>_{L^2(\bbar\lambda)}\\
&=&\frac12\Big[\<E^X_t,\bbar{\Phi}'_{T-t} E^X_t\>_{L^2(\bbar\lambda)}+\<E^X_t,\bbar{\Phi}_{T-t} \wh E^X_t\>_{L^2(\bbar\lambda)}+\<\bbar{\Phi}_{T-t} E^X_t, \wh E^X_t\>_{L^2(\bbar\lambda)}\Big]\\
&=&\frac12\int\int E^X_t(x)E^X_t(y)\Big(\varphi'({T-t},x,y)+(f(x)+f(y))\varphi({T-t},x,y)\Big)\,\bbar\lambda(dx)\,\bbar\lambda(dy)\\
&=&\frac12\int\int E^X_t(x)E^X_t(y)\big(f(x)+\psi({T-t},x)\big)\big(f(y)+\psi({T-t},y)\big)\,\bbar\lambda(dx)\,\bbar\lambda(dy)\\
&=&\frac12\bigg(\int E^X_t(x)\big(f(x)+\psi({T-t},x)\big)\,\bbar\lambda(dx)\bigg)^2,
\end{eqnarray*}
where we have used the Riccati equation \eqref{InfiniteDimRiccatiEqn} and the notation \eqref{psiDefinitionEq} in the fourth step. This proves the assertion.\qed

\bigskip

It follows from the above that a $[0,T]$-admissible strategy $X^*$ with $X^*_0=x$ satisfies
$$\cC(X^*)=V(T,0,x)\le\cC(X)
$$
for all other $[0,T]$-admissible strategies $X$ with $X_0=x$ if $dC^{X^*}_t/dt$ vanishes for a.e. $t\in[0,T]$. Using \eqref{dCXdtEq}, we write this latter condition as 
\begin{equation}\label{X*Cond2}
0=X^*_t+\int E^{X^*}_t(\rho)\theta(T-t,\rho)\,\lambda(d\rho)\qquad\text{for a.e. $t$},
\end{equation}
where
$$\theta(\tau,\rho)=\frac{\varphi_{0}(\tau)(\rho+\psi(\tau,\rho))}{\psi(\tau,0)}-\varphi(\tau,\rho,0).
$$
Then
\begin{equation}\label{ThetaIntEq}
\int\theta(\tau,\rho)\,\lambda(d\rho)=\frac{\varphi_0(\tau)2\bbar\rho}{\psi(\tau,0)}-1.
\end{equation}
Plugging this and \eqref{dEtEquation} into \eqref{X*Cond2} yields that for a.e. $t$
\begin{equation}\label{X*auxEq}
X^*_t=X^*_0\Big(1-\frac{\psi(T-t,0)}{\varphi_0(T-t)2\bbar\rho}\Big)+\frac{\psi(T-t,0)}{\varphi_0(T-t)2\bbar\rho}\int_0^t\int\rho E^{X^*}_s(\rho)\theta(T-t,\rho)\,\lambda(d\rho)\,ds.
\end{equation}
Thus, the left-continuous function $X^*_t$ coincides with an absolutely continuous function for a.e. $t\in[0,T]$. It follows that these two functions coincide for every $t\in(0,T]$. Thus, $E^{X^*}$ is  continuous on $(0,T]$ by \eqref{dEtEquation}, which in turn implies via \eqref{X*auxEq} that  $X^*$ is continuously differentiable throughout $(0,T)$.

When taking the limit $t\downarrow0$ in \eqref{X*auxEq}, we get 
$$X^*_{0+}=x\Big(1-\frac{\psi(T,0)}{\varphi_0(T)2\bbar\rho}\Big),
$$
which gives
\begin{equation}\label{InitialJump}
\Delta X_0^*=-\frac{\psi(T,0)}{2\bbar\rho\varphi_0(T)}\,x.
\end{equation}
That $\Delta X^*_0=\Delta X_T^*$ follows from Remark 2.10 in \citeasnoun{GSS}.

Since  $\psi(t,\rho)$ is continuously differentiable in $t$,  $\theta(t,\rho)$
is also continuously differentiable in $t$. Differentiating \eqref{X*Cond2} with respect to $t\in(0,T)$ yields
\begin{eqnarray*}0&=&\frac{d}{dt}X^*_t+\int \frac{dE^{X^*}_t(\rho)}{dt}\theta(T-t,\rho)\,\lambda(d\rho)-\int E^{X^*}_t(\rho)\theta'(T-t,\rho)\,\lambda(d\rho)\\
&=& \frac{2\bbar\rho\varphi_0(T-t)}{\psi(T-t,0)}\cdot\frac{d}{dt}X^*_t-\int E^{X^*}_t(\rho)\big(\theta'(T-t,\rho)+\rho\theta(T-t,\rho)\big)\,\lambda(d\rho),
\end{eqnarray*}
where we have used \eqref{dEtEquation} and \eqref{ThetaIntEq} in the second
step. This gives
\begin{equation}\label{EtrhoODE}
\frac{d}{dt}X^*_t= \frac{\psi(T-t,0)}{2\bbar\rho\varphi_0(T-t)}\int E^{X^*}_t(x)\big(\theta'(T-t,x)+x\theta(T-t,x)\big)\,\lambda(dx).
\end{equation}

We now want to simplify \eqref{EtrhoODE}. To this end, we use the notation $\bbar\psi(t):=\int x\psi(t,x)\,\lambda(dx).
$ and the formulas
\begin{eqnarray}
\varphi'(t,\rho,0)&=&-\rho\varphi(t,\rho,0)+\frac1{2\bbar\rho}(\rho+\psi(t,\rho))\psi(t,0), \qquad
\varphi_0'(t)=\frac1{2\bbar\rho}\psi(t,0)^2\nonumber\\
\psi'(t,\rho)&=&-\int x^2\varphi(t,x,\rho)\,\lambda(dx)-\rho\psi(t,\rho)+\frac1{2\bbar\rho}(\rho+\psi(t,\rho))(\bbar{\rho^2}+\bbar\psi(t))\label{psiDerivativeEq}\\
\psi'(t,0)&=&-\int x^2\varphi(t,x,0)\,\lambda(dx)+\frac1{2\bbar\rho}\psi(t,0)(\bbar{\rho^2}+\bbar\psi(t)),\nonumber
\end{eqnarray}
Then a tedious computation shows that
\begin{eqnarray*}\theta'(t,\rho)+\rho\theta(t,\rho)=\frac{\varphi_0(t)}{\psi(t,0)}\cdot \Theta(t,\rho),
\end{eqnarray*}
where $\Theta(t,\rho)$ is as in \eqref{ThetaEq}. 
Therefore,
\eqref{EtrhoODE} becomes
\begin{equation}\label{EODE}
\frac{d}{dt}X^*_t=\frac1{2\bbar\rho}\int E^{X^*}_t(x)\Theta(T-t,x)\,\lambda(dx).
\end{equation}
Now we have $E^{X^*}_t=\Delta X^*_0+\int_0^te^{-\rho(t-s)}\frac{d}{dt}X^*_t\,ds$. Plugging this formula back into \eqref{EODE} and using Fubini's theorem yields that $\xi (t)=\frac{d}{dt}X^*_t$ solves the Volterra integral equation \eqref{VolterraEqn}. This is a Volterra integral equation of the second kind with continuous kernel $K(t,s)$ and continuous function $f(t)$. It hence admits a unique continuous solution $x(\cdot)$  \cite[Theorem 3.1]{Linz}. Conversely, given such a solution $x(\cdot)$, we can define a $[0,T]$-admissible strategy $X^*$ via \eqref{InitialJump} and $x(t)=\frac{d}{dt}X^*_t$. Then $X^*$ satisfies  \eqref{X*auxEq} for $t=0+$ as well as \eqref{EtrhoODE} for $t>0$. Integrating \eqref{EtrhoODE}  and reversing the steps made above in deriving \eqref{EtrhoODE} from  \eqref{X*auxEq}  shows that $X^*$ satisfies \eqref{X*auxEq}  for $t\in(0,T]$, and so $X^*$ is optimal. This concludes the proof of Theorem \ref{verificationThm}.\qed

\appendix
\section{Heuristic derivation of the value function}\label{app_heur}

We want to explain here how it is possible to guess the value function
${V}(t,E(\cdot),x)$ introduced in~\eqref{VdefEq}. We start our discussion by
deriving a formula for the costs $\cC(X)$ of a strategy $(X_s,s\in[0,T])$ that 
is arbitrary on $[0,t)$ and optimal on $[t,T]$. We set
\begin{eqnarray*}\wt \cC_t (X)&:=&\cC(X)-\frac12\int_{[0,t)}\int_{[0,t)}G([s-r|)\,dX_r\,dX_s\\
&=&\frac12\int_{[t,T]}\int_{[t,T]}G([s-r|)\,dX_r\,dX_s+\int_{[0,t)}\int_{[t,T]}G(|s-r|)\,dX_s\,dX_r.
\end{eqnarray*}
The rightmost integral can be written as
$$\int_{[t,T]}\int_{[0,t)}G(|s-r|)\,dX_r\,dX_s=\int_{[t,T]}A(s)\,dX_s,
$$
where
$$A(s)=\int E^X_t(\rho)e^{-\rho (s-t)} \lambda(d\rho),\qquad t\le s\le T.
$$
The first-order condition of optimality thus reads
\begin{equation}\label{1stOrderCond}
\int_{[t,T]}G(|s-r|)\, dX_r+A(s)=\nu\qquad\text{for $t\le s\le T$,}
\end{equation}
where $\nu$ is a suitable Lagrange multiplier (compare Theorem 2.11 in \citeasnoun{GSS}).

\begin{lemma}\label{UniquenessLemma} Suppose that $R$ and $\wt R$ are
  functions with finite variation such that $R_T=\wt R_T$ and
$$\int_{[t,T]}G(|s-r|)\,dR_s=\int_{[t,T]}G(|s-r|)\,d\wt R_s\qquad\text{for all $r\in[t,T]$.}
$$
Then $R_s=\wt R_s$ for all $s\in[t,T]$.
\end{lemma}

\Proof We have $\int_{[t,T]}G(|s-r|)\,d(R_s-\wt R_s)=0$ and hence that 
$$\int_{[t,T]}\int_{[t,T]}G(|s-r|)\,d(R_s-\wt R_s)\,d(R_r-\wt R_r)=0,
$$
which implies the assertion in view of the fact that $G$ is strictly positive definite.\qed 

\bigskip
Now suppose that we have auxiliary functions with finite variation $B_t(\rho)$ such that $B_T(\rho)=0$ and
$$\int_{[t,T]}G(|s-r|)\,dB_r(\rho)=e^{-\rho(s-t)} \qquad \text{ for }t\le s\le T.
$$
 We also define 
$$Z_s:=\int E_t^X(\rho)B_s(\rho)\lambda(d \rho) ,$$
so that
\begin{equation}
A(s)=\int_{[t,T]}G(|s-r|)\,dZ_r, \text{ for } t\le s\le T.
\end{equation}
Therefore,
$$\int_{[t,T]}G(|s-r|)\,d(X_r+Z_r)=\nu=\nu \int_{[t,T]}G(|s-r|)\,dB_r(0)\qquad\text{for $t\le s\le T$}.
$$
Lemma~\ref{UniquenessLemma}  hence implies that $X_s+Z_s=\nu B_s(0)\quad\text{for $t\le s\le T$.}$
Hence, we get
$$\nu=\frac{X_t+Z_t}{B_t(0)}, \
X_s=\frac{X_t+Z_t}{B_t(0)}\,B_s(0)-Z_s\qquad\text{for $s\in[t,T]$.}
$$
From these identities we get
\begin{eqnarray*}\wt \cC_t(X)&=&\frac12\int_{[t,T]}\int_{[t,T]}G([s-r|)\,dX_r\,dX_s+\int_{[t,T]}A(s)\,dX_s\\
&=&\frac12 \nu^2\int_{[t,T]}\int_{[t,T]}G(|s-r|)\,dB_r(0)\,dB_s(0)-\nu\int_{[t,T]}\int_{[t,T]}G(|s-r|)\,dB_r(0)\,dZ_s\\
&&+\frac12\int_{[t,T]}\int_{[t,T]}G(|s-r|)\,dZ_r\,dZ_s+\int_{[t,T]}\int_{[t,T]}G(|s-r|)\,dZ_r\,dX_s\\
&=&\frac12\bigg[\frac{(X_t+Z_t)^2}{-B_t(0)}-\int_{[t,T]}\int_{[t,T]}G(|s-r|)\,dZ_r\,dZ_s\bigg],
\end{eqnarray*}
since the first double integral is equal to $-B_t(0)$.
Now we define
$$\varphi(T-t,\rho_1,\rho_2):=\int_{[t,T]}e^{-\rho_1(r-t)}\,dB_r(\rho_2),\qquad
\rho_1,\rho_2\ge 0.
$$
Then, we have
\begin{eqnarray*}
  \varphi(T-t,\rho_1,\rho_2)&=&\int_{[t,T]}\int_{[t,T]}G(|r-s|)\,dB_s(\rho_1)\,dB_r(\rho_2)=\varphi(T-t,\rho_2,\rho_1), \\
\int \varphi(T-t,\rho_1,\rho_2) \, \lambda(d \rho_1) &=&\int_{[t,T]}G(r-t)dB_r(\rho_2)=1.
\end{eqnarray*}
Moreover, we observe that $\varphi(T-t,0,\rho)=-B_t(\rho)$ since
$B_T(\rho)=0$. This
gives in particular that $Z_t=-\int E^X_t(\rho)\varphi(T-t,0,\rho) \lambda(d\rho)$. Besides, we have 
\begin{eqnarray*}\int_{[t,T]}\int_{[t,T]}G(|s-r|)\,dZ_r\,dZ_s&=&\int_{[t,T]}A(s)\,dZ_s\\
&=&\int E^X_t(\rho_1) \int_{[t,T]}e^{-\rho_1(s-t)}\,dZ_s \,\lambda(d\rho_1)\\
&=&\int \int E^X_t(\rho_1)E^X_t(\rho_2)\int_{[t,T]}e^{-\rho_1(s-t)}\,dB_s(\rho_2) \,\lambda(d\rho_1) \,\lambda(d\rho_2)\\
&=&\int \int
E^X_t(\rho_1)E^X_t(\rho_2) \varphi(T-t,\rho_1,\rho_2)
\,\lambda(d\rho_1) \,\lambda(d\rho_2).
\end{eqnarray*}
Thus, we obtain altogether that
$$\wt \cC_t(X)=\frac12\bigg[\frac{\Big(X_t-\int E^X_t(\rho)\varphi(T-t,0,\rho) \lambda(d\rho)\Big)^2}{\varphi(T-t,0,0)}-\int \int
E^X_t(\rho_1)E^X_t(\rho_2) \varphi(T-t,\rho_1,\rho_2)
\,\lambda(d\rho_1) \,\lambda(d\rho_2) \bigg].
$$
This suggests the definition
$$V(t,E(.),X):=\frac12\bigg[\frac{\Big(X-\int E(\rho)\varphi(T-t,0,\rho) \lambda(d\rho)\Big)^2}{\varphi(t,0,0)}-\int \int
E(\rho_1)E(\rho_2) \varphi(t,\rho_1,\rho_2)
\,\lambda(d\rho_1) \,\lambda(d\rho_2)  \bigg],
$$
so that we have
\begin{equation}\label{cost_opt_tT}
  \cC(X)=\int_{[0,t)} \int E^X_s(\rho) \,\lambda(d\rho)\,dX_s+\frac12\sum_{s<t}(\Delta X_s)^2+V(T-t,E_t(.),X_t).
\end{equation}
Thus, we have obtained the formula given by~\eqref{VdefEq}, and it remains to
explain why $\varphi$ should solve a Riccati equation. To do so, we consider
an arbitrary strategy $(X_s,s\in[0,T])$ and consider the cost~\eqref{cost_opt_tT}, which is the cost of the strategy that is equal
to~$X$ on $[0,t)$ and optimal on $[t,T]$. To make the dependence on $t$ explicit, we denote this cost by $\cC_t(X)$. To simplify things, we will focus on the
particular case \eqref{pol_G} of a
discrete measure $\lambda(dx)=\sum_{i=0}^d \lambda_i \delta_{\rho_i}(dx)$, with 
$\rho_0=0<\rho_1<\dots<\rho_d$, $\lambda_i\ge 0$, and
$\sum_{i=0}^d\lambda_i=1$. With this choice, $V$ only depends on $E(\rho_i)$,
$0 \le i\le d$. We introduce the following notations:
\begin{eqnarray*}
  \varphi_{ij}(t)&=&\varphi(t,\rho_i,\rho_j), 0\le i,j\le d,\\
V(t,E_0,\dots,E_d,X)&=&\frac12\bigg[\frac1{\varphi_{00}(t)}\Big(X-\sum_{i=0}^d\lambda_iE_i\varphi_{0i}(t)\Big)^2-\sum_{i=0}^d\sum_{j=0}^d
\lambda_i \lambda_jE_iE_j\varphi_{ij}(t)\bigg], \\
\bbar E^X_t&=&\sum_{i=0}^d \lambda_i
E^X_t(\rho_i),\\
\cC_t(X)&=&\int_{[0,t)} \bbar E^X_s \,dX_s+\frac12\sum_{s<t}(\Delta X_s)^2+V(T-t,E^X_t(\rho_0),\dots,E^X_t(\rho_d)
,X_t).
\end{eqnarray*}

\begin{lemma}We have
$\Delta \cC_t(X)=0$ for all $t \in [0,T]$.
\end{lemma}

\Proof Note that $\Delta X_t=\Delta E^X_t(\rho)$. Hence it is clear that $\Delta \cC_t(X)=0$ if $\Delta X_t=0$. Now suppose that $\Delta X_t\neq0$. Then
\begin{equation}
\Delta \cC_t(X) =\sum_{i=0}^d \lambda_i E^X_t(\rho_i) \Delta
X_t+\frac12(\Delta X_t)^2+\Delta \wt \cC_t(X).
\end{equation}
On the other hand, we have
\begin{eqnarray*}
\lefteqn{V(t,E_0+\delta,\dots,E_d+\delta,X+\delta)}\\&=&\frac12\bigg[\frac1{\varphi_{00}(t)}\Big(X+\delta-\sum_{i=0}^d\lambda_i(E_i+\delta)\varphi_{0i}(t)\Big)^2-\sum_{i=0}^d\sum_{j=0}^d \lambda_i \lambda_j(E_i+\delta)(E_j+\delta)\varphi_{ij}(t)\bigg]\\
&=& V(t,E_0,\dots,E_d,X)-\delta \sum_{j=0}^d\sum_{i=0}^d \lambda_i \lambda_jE_i\varphi_{ij}(t)-\frac{\delta^2}2\sum_{j=0}^d\sum_{i=0}^d \lambda_i \lambda_j\varphi_{ij}(t)\\
&=&V(t,E_0,\dots,E_d,X)-\delta \sum_{i=0}^d \lambda_i E_i -\frac{\delta^2}2.
\end{eqnarray*}
Here we have used the facts that $\sum_i\lambda_i=1$ and $\sum_i\lambda_i\varphi_{ij}=1$. Putting everything together yields the assertion. \qed

\bigskip
We can now focus on infinitesimal variations, and we denote $V_t:=\partial
V/\partial t$, $V_i:=\partial
V/\partial E_i$ and $V_X:=\partial
V/\partial X$. We have, when $\Delta X_t=0$,
\begin{eqnarray*}
d\cC_t(X)&=&\bbar E^X_t \,dX_t-V_t\,dt+\sum_{i=0}^dV_i\,dE^X_t(\rho_i)+V_X\,dX_t\\
&=&\Big(\bbar E^X_t+\sum_{i=0}^dV_i+V_X\Big)\,dX_t-\Big(V_t+\sum_{i=0}^d\rho_iE^X_t(\rho_i)V_i\Big)\,dt\\
\end{eqnarray*}
By simple calculations, we get $\bbar E^X_t+\sum_{i=0}^dV_i+V_X=0$, and our expression simplifies to
$$d\cC_t(X)=-\Big(V_t+\sum_{i=0}^d\rho_iE^X_t(\rho_i)V_i\Big)\,dt.
$$

Let us now calculate $V_t$:
$$V_t=-\frac{\varphi_{00}'}{2\varphi_{00}^2}\Big(X-\sum_{j=0}^d\lambda_jE_j\varphi_{0j}\Big)^2-\sum_{i=0}^d\lambda_iE_i\Big(X-\sum_{j=0}^d\lambda_jE_j\varphi_{0j}\Big)\frac{\varphi'_{0i}}{\varphi_{00}}-\frac12\sum_{i,j=0}^d\lambda_i\lambda_jE_iE_j\varphi_{ij}'.
$$
To simplify computations, we define
\begin{equation}\label{rho-1Eq}
\lambda_{-1}:=1\qquad \text{ and}\qquad  E_{-1}:=\frac{X-\sum_{j=0}^d\lambda_jE_j\varphi_{0j}}{\varphi_{00}}
\end{equation}
as well as $\varphi_{-1i}:=\varphi_{0i}$ and $\varphi_{-1-1}=\varphi_{00}$. Then
$$V_t=-\frac12\sum_{i,j=-1}^d\lambda_i E_i\lambda_j E_j\varphi_{ij}'.
$$
With $\rho_{-1}:=\rho_0=0$, we get
$V_i=-\lambda_i\Big(\lambda_{-1} E_{-1}\varphi_{-1i}+\sum_{j=0}^d\lambda_j  E_j\varphi_{ij}\Big)=-\lambda_i\sum_{j=-1}^d\lambda_j E_j\varphi_{ij}
$ and therefore
$$\sum_{i=0}^d\rho_iE_i V_i=-\sum_{i,j=-1}^d\frac{\rho_i+\rho_j}2\lambda_i E_i\lambda_j E_j\varphi_{ij}.
$$
Altogether, we obtain
\begin{eqnarray*}
d\cC_t(X) =\frac12\Big(\sum_{i,j=-1}^d\lambda_i E^X_t(\rho_i)\lambda_j E^X_t(\rho_j)\big(\varphi_{ij}'+({\rho_i+\rho_j})\varphi_{ij}\big)\Big)\,dt,
\end{eqnarray*}
where $E^X_t(\rho_{-1})$ is defined according to \eqref{rho-1Eq}.

\bigskip

Thus, we arrive at  the following quadratic form $
\frac12\sum_{k,l=0}^dE_kE_l\lambda_k\lambda_l(\varphi'_{kl}+(\rho_k+\rho_l)\varphi_{kl})$. Since
we should have $d\cC_t(X) \ge 0$ with $d\cC_t(X) = 0$ for the optimal
strategy, this quadratic form should be nonnegative with rank one. Indeed,
since the control~$X$ is of dimension one, it would not be possible to make
$d\cC_t(X) = 0$ if the rank of the quadratic form were higher than two. We are
now going to write the conditions that ensures that this quadratic form is
nonnegative with rank one. To do so, 
We introduce the new coordinates $(\Delta_0,\dots, \Delta_d)$ such that
$$E_0=\Delta_0,\quad E_1=\Delta_0+\Delta_1,\quad\cdots\quad E_d=\Delta_0+\Delta_d.
$$
In these coordinates, our quadratic form becomes
\begin{eqnarray*}
\frac12\Delta_0^2\lambda_0^2(\varphi_{00}'+2\rho_0\varphi_{00})+\sum_{l=1}^d\Delta_0(\Delta_0+\Delta_l)\lambda_0\lambda_l(\varphi_{0l}'+(\rho_0+\rho_l)\varphi_{0l})\\+\frac12\sum_{k=1}^d\sum_{l=1}^d(\Delta_0+\Delta_k)(\Delta_0+\Delta_l)\lambda_k\lambda_l(\varphi_{kl}'+(\rho_k+\rho_l)\varphi_{kl})
\end{eqnarray*}
After some calculations, we get that the coefficient for $\Delta_0^2$,
$\Delta_0\Delta_l$, $\Delta_k\Delta_l$ and $\Delta_l^2$ (for $1 \le k,l\le d$)
are respectively $\bbar\rho$,
$\lambda_l\Big(\rho_l+\sum_{k=0}^d\lambda_k\rho_k \varphi_{kl}\Big)$,
$\lambda_k\lambda_l(\varphi_{kl}'+(\rho_k+\rho_l)\varphi_{kl})$ and $\frac12 \lambda_l^2(\varphi_{ll}'+2\rho_l\varphi_{ll})$.

Thus, the matrix ${\bm Q}$ for the quadratic form has coefficients
\begin{eqnarray*}{\bm Q}_{00}&=&\bbar\rho,\\
{\bm Q}_{0l}={\bm Q}_{l0}&=&\frac{\lambda_l}2\Big(\rho_l+\sum_{k=0}^d\lambda_k\rho_k\varphi_{kl}\Big),\\
{\bm Q}_{kl}={\bm Q}_{lk}&=&\frac12\lambda_k\lambda_l(\varphi_{kl}'+(\rho_k+\rho_l)\varphi_{kl})\qquad\text{if $k,l\ge2$, $k\neq l$,}\\
{\bm Q}_{ll}&=&\frac12 \lambda_l^2(\varphi_{ll}'+2\rho_l\varphi_{ll}).
\end{eqnarray*}
Since ${\bm Q}$ is of rank one, the determinant of the matrices
$$\left( \begin{array}{cc}
{\bm Q}_{00} & {\bm Q}_{0l}  \\
{\bm Q}_{l0} & {\bm Q}_{ll}  \end{array} \right), \ \left( \begin{array}{cc}
{\bm Q}_{00} & {\bm Q}_{0l}  \\
{\bm Q}_{k0} & {\bm Q}_{kl}  \end{array} \right)
$$
must vanish for $l=1,\dots, d$ and $k<l$. That gives, respectively,
\begin{eqnarray}
 \varphi_{ll}'+2\rho_l\varphi_{ll}&=&\frac1{2{\bbar\rho}}\Big(\rho_l+\sum_{k=0}^d\lambda_k\rho_k\varphi_{kl}\Big)^2,
 \nonumber \\
\varphi_{kl}'+(\rho_k+\rho_l)\varphi_{kl}&=&\frac1{2{\bbar\rho}}\Big(\rho_l+\sum_{i=0}^d\lambda_i\rho_i\varphi_{il}\Big)\Big(\rho_k+\sum_{j=0}^d\lambda_j\rho_j\varphi_{kj}\Big), \label{derivation_ricc}\end{eqnarray}
which gives precisely the Riccati equation. 
Thus, equation~\eqref{derivation_ricc} holds for $1\le k, l \le d$. In fact, the choice of $E_0=\Delta_0$ is arbitrary. Had we chosen $E_i=\Delta_i$ for some $i>0$, we had obtained~\eqref{derivation_ricc} for $k,l\neq i$. Therefore~\eqref{derivation_ricc} holds in fact for all $k,l=0,\dots,d$.

\parskip-0.5em\renewcommand{\baselinestretch}{0.9}\normalsize
\bibliographystyle{agsm}
\bibliography{MarketImpact}
\end{document}